\documentclass[11pt]{amsart}
\usepackage{amsfonts}
\usepackage{amssymb}
\usepackage[arrow,matrix]{xy}
\usepackage{amsmath,amssymb, bbm, amscd, amsthm,mathrsfs,hyperref,mathtools}
\usepackage[titletoc]{appendix}
\usepackage{tikz-cd}

\DeclareMathAlphabet{\mathpzc}{OT1}{pzc}{m}{it}

\theoremstyle{plain}
\newtheorem{thm}{Theorem}[section]

\newtheorem{lem}[thm]{Lemma}
\newtheorem{prop}[thm]{Proposition}
\theoremstyle{definition}
\newtheorem{defn}[thm]{Definition}

\newtheorem{lem-defn}[thm]{Lemma-Definition}
\newtheorem{thm0}{Theorem}
\newtheorem{question0}[thm0]{Question}

\def\dt{\delta}
\def\gm{\gamma}

\def\vph{\varphi}
\def\lmd{\lambda}

\def\vps{\varepsilon}

\def\bdt{\Delta}

\def\mc{\mathcal}

\def\t{\text}
\def\it{\textit}

\def\ot{\otimes}
\def\op{\oplus}
\def\se{\leqslant}
\def\le{\geqslant}
\def\bxt{\boxtimes}
\def\sq{\square}

\def\ra{\rightarrow}
\def\lra{\longrightarrow}
\def\la{\leftarrow}
\def\xra{\xrightarrow}

\def\Hom{\operatorname {Hom}}
\def\Ext{\operatorname {Ext}}

\def\dim{\operatorname {dim}}
\def\id{\operatorname {id}}

\def\H{\operatorname {H}}

\def\injdim{\operatorname{injdim}}

\def\GL{\operatorname {GL}}
\def\SL{\operatorname {SL}}
\def\tr{\operatorname {tr}}
\def\ob{\operatorname {ob}}

\def\cd{\operatorname {cd}}

\def\kk{\mathbbm{k}}

\def\NN{\mathbb{N}}

\def\dd{\mathbbm{D}}

\begin{document}
\title[Calabi-Yau property]{\bf  Calabi-Yau property of Quantum Groups of $\GL(2)$ Representation Type}

\author{Xiaolan YU}
\address {Xiaolan YU\newline Department of Mathematics, Hangzhou Normal University, Hangzhou, Zhejiang 310036, China}
\email{xlyu@hznu.edu.cn}

\author{Xingting WANG}
\address{Xingting Wang\newline Department of Mathematics, Louisiana State University, Baton Rouge, Louisiana 70803, USA} 

\email{xingtingwang@lsu.edu}

\date{}

\begin{abstract}
In this paper, we prove that the quantum groups $\mc{G}(A,B)$ introduced by Mrozinski in \cite{mro} and their Hopf-Galois objects are twisted Calabi-Yau algebras, and give their Nakayama automorphisms explicitly. 
\end{abstract}



\maketitle

\section*{Introduction}\label{0}

In order to classify the cosemisimple Hopf algebras whose corepresentation semi-ring is isomorphic to that of $
\GL(2)$, C. Mrozinski introduced  a new family of Hopf algebras \cite{mro} as follows. 

Throughout, let  $\kk$ be an algebraically closed field.  Let $n\in \mathbb{N} $ with $n\le 2$ and let $A,B\in \GL_n(\kk)$. The algebra $\mathcal{G}(A,B)$ is presented by generators $(u_{ij})_{1 \se i,j\se n}$, $\dd, \dd^{-1}$ subject to relations
\begin{equation}\label{relations}u^tAu=A\dd,\;\;\;uBu^t=B\dd,\;\;\; \dd \dd^{-1}=1=\dd^{-1}\dd,\end{equation}
where $u$ is the matrix $(u_{ij})_{1\se i,j\se n}$ and $u^t$ denotes its transpose. This algebra has a natural Hopf algebra structure described in \cite[Proposition 2.1]{mro}. Let $$A_q=
\begin{pmatrix}
	0  &  1\\
	-q  &  0\\
\end{pmatrix},
\ \text{for some}\ q\in \kk^\times.$$
The Hopf algebra $\mc{G}(A_q,A_q^{-1})$ is just $\mc{O}(\GL_q(2))$,  namely, the function algebra on the quantum group $\GL_{q}(2)$.  For a matrix $A\in \GL_n(\kk)$,  in \cite{dvl}, Dubois-Violette and Launer defined the Hopf algebra $\mc{B}(A)$ as the quantum automorphism group of the non-degenerate bilinear form associated to $A$. There is a surjective Hopf algebra morphism from $\mc{G}(A,A^{-1})$ to the Hopf algebra $\mc{B}(A)$. In view of the definition, the Hopf algebras $\mathcal{G}(A,B)$ might be viewed as a generalization of 
the Hopf algebras $\mc{B}(A)$.

In \cite[Theorem 6.1 and Corollary 6.3]{bi}, Bichon proved that the Hopf algebras $\mc{B}(A)$ are smooth of dimension 3 and satisfy the Poincar\'{e} duality. It is equivalent to say that they are twisted Calabi-Yau (CY for short) algebras  of dimension 3 (Definition \ref{defn tcy}).  CY algebras were introduced by Ginzburg in 2006 \cite{g2} to transport the geometry of CY manifolds to noncommutative algebraic geometry. Nowadays, CY algebras and their categorifications are deeply entangled with noncommutative geometry, cluster algebras, topological field theories and representation theory of noncommutative algebras. Twisted (or ``skew'') CY algebras are a generalization of CY algebras.  Associated to a twisted CY algebra, there exists a so-called Nakayama automorphism. This algebra automorphism is unique up to an inner one. A CY algebra is a  twisted CY algebras whose Nakayama automorphism is inner.   Actually, the twisted CY property has been studied for many years, even before the definition of CY algebras, under the notion of Artin-Schelter (AS for short) regularity. AS-regular algebras play a central role in the study of noncommutative projective geometry. A connected graded algebra is graded twisted CY if and only if it is AS-regular (cf. \cite[Lemma 1.2]{rrz}). Therefore, twisted CY algebras can also be viewed as a generalization of AS-regular algebras. In particular in the case of Hopf algebras, it is proved in \cite{bz} (cf. \cite[Lemma 1.3]{rrz}) that Noetherian AS-regular Hopf algebras are twisted CY algebras (the twisted CY property is called ``rigid Gorenstein'' in \cite{bz}).

The main aim of this paper is to prove that the algebras $\mc{G}(A,B)$ are twisted CY algebras. We construct a free Yetter-Drinfeld resolution (Definition \ref{defn ydres}) of the trivial Yetter-Drinfeld module $\kk$ over $\mc{G}(A,B)$ (Theorem \ref{main thm}). As a consequence,  we obtain the following  main theorem.  
\begin{thm0}\label{thm0}
	Let $A,B\in \GL_n(\kk)\;(n\le 2)$ such that $B^tA^tBA=\lambda I_n$ for some $\lmd\in \kk^\times$. The algebra $\mc{G}(A,B)$ is a twisted CY algebra with Nakayama automorphism $\mu$ defined by $\mu(u)= (A^t)^{-1}A u B^tB^{-1}$ and $\mu(\dd^{\pm 1})=\dd^{\pm 1}$.
\end{thm0}

In the study of homological property, a general classical problem is to find  homological invariants under various kind of ``deformations''. Motivated by the discussion of the twisted CY property of the Hopf algebras $\mc{B}(A)$, in  \cite{wyz}, the authors asked the following question:
\begin{question0}\label{Q1}
	Let $H$ and $L$ be two Hopf algebras such that their  comodule categories are monoidally equivalent. Suppose $H$ is twisted CY. Is $L$ always twisted CY?
\end{question0}
Let $m,n\in \NN , n,m\le 2$ and let $A,B\in \GL_n(\kk),C,D\in \GL_m(\kk)$ such that $B^tA^tBA=\lmd I_n$, $D^tC^tDC=\lmd I_m$ and $\tr (AB^t)=\tr(CD^t)=\mu$ for some $\lmd, \mu\in \kk^{\times}$. Then the comodule categories over  $\mc{G}(A,B)$ and $\mc{G}(C,D)$ are monoidally equivalent \cite{mro}. Theorem \ref{thm0} shows that the algebras $\mc{G}(A,B)$ provide an interesting example to study the above question.

As a corollary of Theorem \ref{thm0}, we prove that  Hopf-Galois objects of the algebras $\mc{G}(A,B)$  are also twisted CY algebras (Theorem \ref{thm galois}). 

Another application of the free Yetter-Drinfeld resolution we construct concerns Gerstenhaber-Schack cohomology (\cite{gs1,gs2}), a cohomology theory specific to Hopf algebras. It has been useful in proving some fundamental results in Hopf algebra theory.  Bialgebra cohomology is a special instance of Gerstenhaber-Schack cohomology. When the algebra $\mc{G}(A,B)$ is cosemisimple, we compute its bialgebra cohomology and we prove that the  Gerstenhaber-Schack cohomological dimension of $\mc{G}(A,B)$ is 4. 

\section{Preliminaries}\label{1}
Throughout this paper, we work over an algebraically closed fixed field $\kk$. Unless stated otherwise all algebras and vector spaces are over $\kk$. The unadorned tensor $\otimes$ means $\otimes_\kk$ and $\Hom$ means $\Hom_\kk$. For an algebra $A$, we denote by $A^{op}$ the opposite algebra of $A$ and $A^e$ the enveloping algebra $A\ot A^{op}$ of $A$. An $A$-bimodule can be identified with an $A^e$-module.

For an $A$-bimodule $M$ and an algebra automorphism $\mu$, we write $M^\mu$ for the  $A$-bimodule defined by $M^\mu\cong M$ as vector spaces, and 
$a\cdot m \cdot b=am\mu(b)$, for all $a,b\in A$ and $m\in M$. Similarly, we have the notion of ${}^\mu M$. It is well-known that   $A^\mu\cong A$ as $A$-bimodules if and only
if $\mu$ is an inner automorphism of $A$.

For a Hopf algebra, we use  Sweedler's (sumless) notation for its comultiplication and its coactions. The category of right $H$-comodule  is denoted by $\mc{M}^H$. Let $H$ be a Hopf algebra, and $\xi:H\ra\kk$ an algebra
homomorphism. We use $[\xi]^l$ to denote the left winding automorphism of
$\xi$ defined by
$$[\xi]^l(h)=\xi(h_{(1)})h_{(2)},$$ for any $h\in H$. The right winding automorphism $[\xi]^r$ of $\xi$ can be defined similarly. That is,   for any  $h\in H$,
$$[\xi]^r(h)=\xi(h_{(2)})h_{(1)}.$$
It is easy to check that both left and right winding automorphisms are algebra automorphisms of $H$ (cf. \cite[\S 2.5]{bz}).

\subsection{Cogroupoid} The monoidally equivalence between comodule categories of Hopf algebras can be described by the language of cogroupoids. 
\begin{defn}A \it{cogroupoid} $\mc{C}$ consists of:
	\begin{itemize}
		\item A set of objects $\ob(\mc{C})$.
		\item For any $X,Y\in \ob(\mc{C})$, an algebra $\mc{C}(X,Y)$.
		\item For any $X,Y,Z\in \ob(\mc{C})$, algebra homomorphisms
		$$\bdt^Z_{XY}:\mc{C}(X,Y)\ra \mc{C}(X,Z)\ot \mc{C}(Z,Y) \t{ and } \vps_X:\mc{C}(X,X)\ra \kk$$
		and linear maps
		$$S_{X,Y}:\mc{C}(X,Y)\longrightarrow \mc{C}(Y,X)$$
		such that for any $X,Y,Z,T\in \ob(\mc{C})$, the following diagrams commute:
\[
\xymatrix{
\mc{C}(X,Y)\ar^-{\bdt^Z_{X,Y}}[rr]\ar^{\bdt^T_{X,Y}}[d] && \mc{C}(X,Z)\ot \mc{C}(Z,Y)\ar^-{\bdt^T_{X,Z}\ot \id}[d]\\
\mc{C}(X,T)\ot \mc{C}(T,Y)\ar^-{\id\ot \bdt^Z_{T,Y}}[rr]&&\mc{C}(X,T)\ot \mc{C}(T,Z)\ot \mc{C}(Z,Y), }
\]
\xymatrix
{\mc{C}(X,Y)\ar@{=}[rd]\ar^{\bdt^Y_{X,Y}}[d]&\\
\mc{C}(X,Y)\ot\mc{C}(Y,Y)\ar^-{\id \ot \vps_Y}[r]&\mc{C}(X,Y),}\hspace{4mm}
\xymatrix
{\mc{C}(X,Y)\ar@{=}[rd]\ar^{\bdt^X_{X,Y}}[d]&\\
\mc{C}(X,X)\ot\mc{C}(X,Y)\ar^-{\vps_X \ot \id }[r]&\mc{C}(X,Y).} 
	\end{itemize}
	and
	\small \[\xymatrix{\mc{C}(Y,X) & \kk \ar[l]_-u & \mc{C}(X,X)\ar[d]_{\bdt_{X,X}^Y}\ar[r]^-{\vps_X}\ar[l]_-{\epsilon_X} &\kk\ar[r]^-u&\mc{C}(X,Y)\\
\mc{C}(Y,X)\ot\mc{C}(Y,X) \ar[u]^m &&\mc{C}(X,Y)\ot\mc{C}(Y,X)\ar[rr]^{\id\ot S_{Y,X}} \ar[ll]_-{S_{X,Y} \otimes \id} &&\mc{C}(X,Y)\ot\mc{C}(X,Y). \ar[u]^m}\]
\end{defn}
The axioms of a cogroupoid are dual to the axioms defining a groupoid. Note that $\mc{C}(X,X)$ is a Hopf algebra for any object $X\in \mc{C}$. A cogroupoid $\mc{C}$ is said to be \it{connected} if $\mc{C}(X,Y)$ is a nonzero algebra for any $X,Y\in \ob(\mc{C})$.

We use Sweedler's notation for cogroupoids. Let $\mc{C}$ be a cogroupoid. For any $a^{X,Y}\in \mc{C}(X,Y)$, we write
$$\bdt^Z_{X,Y}(a^{X,Y})=a^{X,Z}_{(1)}\ot a^{Z,Y}_{(2)}.$$

The following proposition describes properties of the ``antipodes." For other properties of cogroupoids, we refer the reader to \cite{bi1}.
\begin{prop}\cite[Proposition 2.13]{bi1}
Let $\mc{C}$ be a cogroupoid and let $X,Y\in\ob(\mc{C})$.
\begin{enumerate}
\item[(i)] $S_{Y,X}:\mc{C}(Y,X)\ra \mc{C}(X,Y)^{op}$ is an algebra homomorphism.
\item[(ii)] For any $Z\in \ob(\mc{C})$ and $a^{Y,X}\in\mc{C}(Y,X)$,
$$\bdt_{X,Y}^Z(S_{Y,X}(a^{Y,X}))=\sum S_{Z,X}(a_2^{Z,X})\ot S_{Y,Z}(a_1^{Y,Z}).$$
\end{enumerate}
\end{prop}

For a detailed introduction about cogroupoids, we refer to \cite{bi1}.

In \cite{sch}, Schauenburg showed that the categories of comodules over two Hopf algebras $H$ and $L$ are monoidally equivalent, called Morita-Takeuchi equivalent, if and only if there exists an $H$-$L$-bi-Galois object between them. Later, Bichon \cite{bi1} introduced the notion of a cogroupoid to provide a categorial context for Hopf-(bi)Galois objects.

\begin{thm}\cite[Theorem 2.10, 2.12]{bi1}\label{thm ga}
	Let $\mc{C}$ be a connected cogroupoid. Then for any $X,Y\in\mc{C}$, we have equivalences of monoidal categories that are inverse of each other
	$$\begin{array}{rclcrcl}\mc{M}^{\mc{C}(X,X)}&\cong^\ot& \mc{M}^{\mc{C}(Y,Y)}&&\mc{M}^{\mc{C}(Y,Y)}&\cong^\ot& \mc{M}^{\mc{C}(X,X)}\\
		V&\longmapsto&V\sq_{\mc{C}(X,X)}\mc{C}(X,Y)&&V&\longmapsto&V\sq_{\mc{C}(Y,Y)}\mc{C}(Y,X)
	\end{array}$$
	Conversely, if $H$ and $L$ are Hopf algebras such that $\mc{M}^H\cong ^\ot \mc{M}^L$, then there exists a connected cogroupoid with 2 objects $X,Y$ such that $H=\mc{C}(X,X)$ and $L=\mc{C}(Y,Y)$.
\end{thm}


\subsection{Calabi-Yau algebras}
In this subsection, we recall the definition of (twisted) Calabi-Yau algebras.
\begin{defn}\label{defn tcy}An algebra $A$ is called  a \it{twisted Calabi-Yau algebra  of dimension
		$d$} if
	\begin{enumerate}\item[(i)] $A$ is \it{homologically smooth}, that is, as an $A^e$-module, $A$ has
		a finitely generated projective resolution of finite length; 
		\item[(ii)] There is an automorphism $\mu$ of $A$ such that
		\begin{equation}\Ext_{A^e}^i(A,A^e)\cong\begin{cases}0,& i\neq d;
				\\A^\mu,&i=d,\end{cases}\end{equation}
		as $A^e$-modules.
	\end{enumerate}
	If such an automorphism $\mu$ exists, it is unique up to an inner automorphism and is called a \it{Nakayama automorphism} of $A$. When a Nakayama automorphism of a twisted Calabi-Yau algebra is an inner automorphism, it is called a \it{Calabi-Yau algebra}. In what follows, Calabi-Yau is abbreviated to CY for short. 
\end{defn}

As mentioned in the introduction, twisted CY algebras are closely related to Artin-Schelter (AS) regular algebras.

\begin{defn}\label{defn as} Let $H$ be a Hopf algebra, and let $\kk$ also denote the trivial $H$-module. 
	\begin{enumerate}
		\item[(i)] The Hopf algebra $H$ is said to be \it{left
			AS-Gorenstein}, if
		\begin{enumerate}
			\item $\injdim  {_HH}=d<\infty$,
			\item $\dim\Ext_H^i({_H\kk},{_HH})=\begin{cases}0, &i\neq d;\\1&i= d.\end{cases}$
		\end{enumerate}
		A \it{right AS-Gorenstein} Hopf algebra can be defined similarly. 
		\item[(ii)] If $H$ is both left and right AS-Gorenstein, then $H$ is called \textit{AS-Gorenstein}.
		\item[(iii)] If, in addition, the global dimension of $H$ is finite, then $H$ is called  \it{AS-regular}.
	\end{enumerate}
\end{defn}

Compare with \cite[Definition 1.2]{bz}, we do not require the Hopf algebra $H$ to be Noetherian.  For AS-regularity,  the right global dimension always equals the left global dimension. When $H$ is AS-Gorenstein and homologically smooth, the right injective dimension always equals the left injective dimension, which are both given by the integer $d$ such that $\Ext_{H^e}^d(H,H^e)\neq 0$.
we refer to \cite[Remark 2.1.5]{wyz} for an explanation. 

Let $H$ be a left AS-Gorenstein Hopf algebra of injective dimension $d$.  The one dimensional right $H$-module $\Ext^d_H({}_H\kk, {_HH})$ is usually denoted by $\int^l_H$, and is called the \emph{left  homological integral} of $H$. Let $x$ be a nonzero element in $\int^l_H$,  the right $H$-action defines an algebra homomorphism $\xi:H\ra \kk$ by $x\cdot h=\xi(h)x$, for any $h\in H$. That is,  $\int^l_H\cong \kk_{\xi}$ as right $H$-modules.

If $H$ is a right AS-Gorenstein Hopf algebra of injective dimension $d$, one may define the \emph{right homological integral} $\int^r_H=\Ext^d_H(\kk_H, {H_H})$ similarly. 
There is an algebra homomorphism $\eta:H\ra \kk$, such that $\int^r_H\cong {_\eta\kk}$ as left $H$-modules.

The following lemma summarizes the relationship between twisted CY Hopf algebras and AS-Gorenstein  Hopf algebras (not necessarily Noetherian). They are probably well-known. We list them here for the sake of readers' convenience. 
\begin{lem}\label{lem AS}
	Let $H$ be a Hopf algebra with bijective antipode. Then the followings are equivalent, 
	\begin{enumerate}
		\item[(i)]  $H$ is twisted CY;
		\item[(ii)]  $H$ is left AS-Gorenstein and the left trivial module $\kk$ admits a bounded projective resolution with each term finitely generated;
		\item[(iii)]  $H$ is right AS-Gorenstein and the right trivial module $\kk$ admits a bounded projective resolution with each term finitely generated.
	\end{enumerate}
		Now suppose $H$ is twisted CY. 
		\begin{enumerate}
		\item[(iv)] Let  $\mu$ be a Nakayama automorphism of $H$. We have $\int^l_H\cong \kk_\xi$ as right $H$-modules and $\int^r_H\cong {}_\eta\kk$ as left $H$-modules, where $\xi$ and $\eta$ are algebra homomorphism defined by $\xi=\vps \mu$ and $\eta=\vps \mu^{-1}$, respectively. 
		\item[(v)] If we know that  $\int^l_H\cong \kk_\xi$ as right $H$-modules, where $\xi:H\ra \kk$ is an algebra homomorphism, a Nakayama automorphism of $H$ is given by $\mu=S^2[\xi]^l$. Alternatively,  the algebra automorphism $S^{-2}[\eta S]^r$ is also a Nakayama automorphism $H$, where $\eta:H\ra \kk$  is the algebra homomorphism such that $\int^r_H\cong {_\eta\kk}$ as left $H$-modules.
			\end{enumerate}
			
\end{lem}

\proof The equivalence of (i), (ii) and (iii) is  \cite[Corollary 2.1.7]{wyz}. Therefore, if  $H$ is a twisted CY algebra, then $H$ is both left and right AS-Gorenstein, that is, $H$ is AS-Gorenstein. 

(iv) follows from Lemma 2.15 and Remark 2.16 in \cite{yvz}.

(v) Let $\xi$ be the algebra homomorphism such that  $\int^l_H\cong \kk_{\xi}$ as right $H$-modules. From the proof Proposition 2.1.6 in \cite{wyz}, we can see that 
$$\Ext_{H^e}^i(H,H^e)\cong\begin{cases}0,& i\neq d;
		\\H^\mu,&i=d,\end{cases}$$
where the automorphism $\mu$ is defined by $\mu(h)=\xi(h_{(1)})S^2(h_{(2)})$. That is, $S^2[\xi]^l$ is a Nakayama automorphism of $H$. Similarly,  let $\eta:H\ra \kk$  be the algebra homomorphism such that $\int^r_H\cong {_\eta\kk}$ as left $H$-modules, we have that  $$\Ext_{H^e}^i(H,H^e)\cong\begin{cases}0,& i\neq d;
	\\{}^{^{{S^{2}[\eta]^r}}}H,&i=d.\end{cases}$$
Therefore, $({S^{2}[\eta]^r})^{-1}={S^{-2}[\eta S]^r}$ is also a Nakayama automorphism of $H$. \qed

\subsection{Yetter-Drinfeld modules}

	Let $H$ be a Hopf algebra. A
	\textit{(right-right) Yetter-Drinfeld module} $V$ over $H$ is simultaneously a
	right $H$-module and a right $H$-comodule satisfying the compatibility
	condition $$\dt( v\cdot h) = v_{(0)}\cdot h_{(2)}
	\ot S(h_{(1)}) v_{(1)}h_{(3)},
	$$ for any $v\in V$, $h\in H$. The category of Yetter-Drinfeld  modules is usually denoted by $\mc{YD}^H_H$. 

Let $H$ be a Hopf algebra, and $V$ a right $H$-comodule. The Yetter-Drinfeld module $V\boxtimes H$ is defined by Bichon in \cite[Proposition 3.1]{bi} as follows. The right $H$-module structure is defined
by multiplication on the right side, and the right comodule structure is defined by the linear map
$$\begin{array}{rcl}V\ot H&\ra& V\ot H\ot H\\
	v\ot h&\mapsto &v_{(0)}\ot h_{(2)}\ot S(h_{(1)})v_{(1)} h_{(3)}.\end{array}$$

A  Yetter-Drinfeld module over $H$ is said to be \textit{free} if it is isomorphic to $V\bxt H$ for some right $H$-comodule $V$.

A free Yetter-Drinfeld module is obviously free as a right $H$-module. We call a free Yetter-Drinfeld module $V\bxt H$ \it{finitely generated} if $V$ is finite dimensional.

\begin{defn}\label{defn ydres}
	Let $H$ be a Hopf algebra and let $M\in \mc{YD}^H_H$. A \textit{free Yetter-Drinfeld module resolution} of $M$ is a complex of free  Yetter-Drinfeld modules
	$$\textbf{P}_\bullet:\cdots\ra P_{i+1}\ra P_i\ra \cdots\ra P_1\ra P_0\ra 0$$
	together with  a Yetter-Drinfeld module morphism $\varepsilon:P_0\ra M$ such that the complex
	$$\cdots\ra P_{i+1}\ra P_i\ra \cdots\ra P_1\ra P_0\xra{\epsilon}M\ra 0$$
	is an exact sequence. 
\end{defn}

The monoidal equivalence in Theorem \ref{thm ga} can be extended to categories of Yetter-Drinfeld modules.

	Let $\mc{C}$ be a cogroupoid, $X,Y\in\ob(\mc{C})$ and $V$ a right $\mc{C}(X,X)$-module.  Following from \cite[Proposition 6.2]{bi1}, we know that 
	 $V\ot \mc{C}(X,Y)$ has a right $\mc{C}(Y,Y)$-module structure defined by
		$$ (v\ot a^{X,Y})\la b^{Y,Y} =  v \cdot b_{(2)}^{X,X}\ot S_{Y,X}(b_{(1)}^{Y,X}) a^{X,Y} b_{(3)}^{X,Y}.$$
		Together with the right $\mc{C}(Y,Y)$-comodule structure defined by $1\ot \bdt_{X,Y}^Y$,  $V\ot \mc{C}(X,Y)$ is a Yetter-Drinfeld module over $\mc{C}(Y,Y)$. 
		If moreover $V$ is a Yetter-Drinfeld module, then $V\square_{\mc{C}(X,X)}\mc{C}(X,Y)$ is a  Yetter-Drinfeld submodule of $V\ot \mc{C}(X,Y)$.

\begin{thm}\cite[Theorem 4.4]{bi}\label{lem 4.4}
	Let $\mc{C}$ be a connected cogroupoid. Then for any $X,Y\in \ob(\mc{C})$, the functor 
	$$\begin{array}{rcl}
		\mc{YD}^{\mc{C}(X,X)}_{\mc{C}(X,X)}&\longrightarrow &\mc{YD}^{\mc{C}(Y,Y)}_{\mc{C}(Y,Y)}\\
		V&\longmapsto &V\sq_{\mc{C}(X,X)}\mc{C}(X,Y)
	\end{array}$$
	is an equivalence of monoidal categories. Moreover, we have natural isomorphisms
	$$\begin{array}{rcl}
		(V\sq_{\mc{C}(X,X)}\mc{C}(X,Y))\boxtimes \mc{C}(Y,Y)&\longrightarrow &(V\boxtimes \mc{C}(X,X))\sq_{\mc{C}(X,X)}\mc{C}(X,Y)\\
		v\ot a^{XY} \ot b^{YY}&\longmapsto &v\ot b_{(2)}^{XX}\ot S_{Y,X}(b_{(1)}^{YX})a^{XY}b_{(3)}^{XY}.
	\end{array}$$
	for any $v\in V, a^{XY}\in \mc{C}(X,Y), b^{YY}\in \mc{C}(Y,Y)$.
\end{thm}


\section{A resolution of the trivial Yetter-Drinfeld module}\label{2}
Let $n\in \mathbb N_{\ge 2}$, and $A,B\in \GL_{n}(\kk)$ satisfying $B^tA^tBA=\lambda I_n$  for some nonzero scalar $\lambda$. One can check directly that $A^tB^tAB=\lambda I_n$ as well. Recall the algebra $\mc{G}(A,B)$ that has been defined  in the introduction.  Note that by the relations (\ref{relations}), the following commutation relations hold in $\mc{G}(A,B)$,
$$BA\dd u~=~u\dd BA.$$
That is, $\dd u\dd^{-1}=A^{-1}B^{-1}uBA$.

The algebra $\mc{G}(A,B)$ admits a Hopf algebra structure, with comultiplication defined by
$$\bdt(u_{ij})=\sum_{k=1}^{n} u_{ik}\ot u_{kj} (1\se i,j\se n),\quad \bdt(\dd^{\pm 1})=\dd^{\pm 1}\ot \dd^{\pm 1},$$
with counit defined by 
$$\vps(u_{ij})=\dt_{ij}\ (1\se i,j\se n),\quad \vps(\dd^{\pm 1})=1,$$
and with antipode defined by
$$S(u)=\dd^{-1}A^{-1}u^tA,\quad S(\dd^{\pm 1})=\dd^{\mp 1}.$$ 
The antipode is clearly bijective. Moreover,  it is easy to check that $S^2=\Phi*{\rm id}*\Phi$ with sovereign character $\Phi: \mc{G}(A,B)\to \kk$ defined by the algebra map $u\mapsto BA^t$ and $\dd^{\pm 1}\mapsto \lambda^{\pm 1}$ or 
\[S^2(u)=\dd^{-1}A^{-1}A^tu(A^t)^{-1}A\dd=BA^tu(A^t)^{-1}B^{-1},\quad S^2(\dd^{\pm 1})=\dd^{\pm 1}.\] 

Assume that the characteristic of the field $\kk$ is $0$. Then the Hopf algebras whose corepresentation semi-ring is isomorphic to that of $\GL_2(\kk)$ are exactly the Hopf algebras $\mc{G}(A,B)$ with $A,B\in\GL_n(\kk)$ $(n\le 2)$ satisfying $B^tA^tBA=\lmd I_n$ for some $\lmd\in \kk^{\times}$ and such that any solution of the equation $X^2-\sqrt{\lmd^{-1}}\tr(AB^t)X+1=0$ is generic (\cite[Theorem 1.2]{mro}). 

For the matrix $A_q\in \GL_2(\kk)$ in the introduction, we have $\mc{G}(A_q,A_q^{-1})=\mc{O}(\GL_q(2))$. Let $A,B\in \GL_n(\kk)\;(n\le 2)$ such that $B^tA^tBA=\lambda I_n$ for some $\lmd\in \kk^\times$ and let $q\in\kk^\times$ such that $q^2-\sqrt{\lmd^{-1}}\tr(AB^t)q+1=0$. It is shown in \cite{mro} that the monoidal categories of comodules over $\mc{G}(A,B)$ and $\mc{O}(\GL_q(2))$ are equivalent.

In this section, we will construct a free Yetter-Drinfeld resolution of trivial Yetter-Drinfeld module $\kk$ over $\mc{G}(A,B)$. We first recall the cogroupoid defined in \cite{mro}.

\subsection{Cogroupoid} Let $m,n\in \NN_{\ge 2}$ and let $A,B\in \GL_n(\kk),C,D\in \GL_m(\kk)$. The algebra $\mc{G}(A,B|C,D)$ is defined as follows. It is generated by 
$$(u^{AB,CD}_{ij})_{1 \se i,j\se n},\quad \dd_{AB,CD},\quad \dd_{AB,CD}^{-1}$$ subject to 
\begin{equation}\label{alg}u^tAu=C \dd,\;\;\;uDu^t=B\dd,\;\;\; \dd_{AB,CD} \dd_{AB,CD}^{-1}=1=\dd_{AB,CD}^{-1}\dd_{AB,CD},\end{equation}
where $u$ is the matrix $(u^{AB,CD}_{ij})_{1\se i,j\se n}$. Similarly, in $\mc{G}(A,B|C,D)$ we have 
$$BA\dd u~=~u\dd DC\quad \text{or}\quad CD\dd u^t=u^t\dd AB.$$
That is, $\dd u\dd^{-1}=A^{-1}B^{-1}uDC$ or $\dd u^t\dd^{-1}=D^{-1}C^{-1}u^tAB$. When there is no confusion, we will simply denote $u^{AB,CD}_{ij}$, $\dd_{AB,CD}$ and $\dd_{AB,CD}^{-1}$ by $u_{ij}$, $\dd$ and $\dd^{-1}$, respectively. It is clear that $\mc{G}(A,B|A,B)=\mc{G}(A,B)$.

For any $A,B\in \GL_n(\kk)$, $C,D\in \GL_m(\kk)$ and $X,Y\in \GL_p(\kk)$, define the following maps:
\begin{equation}\label{bdt}\begin{array}{rcl}
		\bdt^{XY}_{AB,CD}:\mc{G}(A,B|C,D)&\longrightarrow &\mc{G}(A,B|X,Y)\ot \mc{G}(X,Y|C,D)\\
		u_{ij}&\longmapsto &\sum_{k=1}^p u_{ik}\ot u_{kj}\\
		\dd^{\pm1}&\longmapsto &\dd^{\pm1}\ot\dd^{\pm1},
\end{array}\end{equation}
\begin{equation}\label{bvps}\begin{array}{rcl}
		\vps_{AB}:\mc{G}(A,B)&\longrightarrow &\kk\\
		u_{ij}&\longmapsto &\dt_{ij}\\
		\dd^{\pm1}&\longmapsto &1,
\end{array}\end{equation}
\begin{equation}\label{bss}\begin{array}{rcl}
		S_{AB,CD}:\mc{G}(A,B|C,D)&\longrightarrow &\mc{G}(C,D|A,B)^{op}\\
		u&\longmapsto &\dd^{-1}A^{-1}u^tC\\
		\dd^{\pm 1}&\longmapsto&\dd ^{\mp 1} 
\end{array}\end{equation}
It is clear that $S_{AB,CD}$ is bijective.

The cogroupoid associated to $\mc{G}(A|B)$ was introduced in \cite[Definition 3.3]{mro} based on \cite[Lemma 3.2]{mro}.  

\begin{defn}
	The cogroupoid $\mc{G}$ is the cogroupoid defined as follows:
	\begin{enumerate}
		\item[(i)] $\ob(\mc{G})=\{(A,B)\in \GL_n(\kk)\times \GL_n(\kk),n\le1\}$.
		\item[(ii)] For $(A,B),(C,D)\in \ob(\mc{G})$, the algebra $\mc{G}(A,B,|C,D)$ is the algebra defined as in (\ref{alg}).
		\item[(iii)] The structural maps $\bdt^\bullet_{\bullet,\bullet}$, $\vps_\bullet$ and $S_{\bullet,\bullet}$ are defined in (\ref{bdt}), (\ref{bvps}) and (\ref{bss}), respectively.
	\end{enumerate}
\end{defn}

\begin{lem}\cite[Corollary 3.6]{mro}\label{lem G E}
	Let $\lmd,\mu\in \kk^\times$. Consider the full subcogroupoid $\mc{G}^{\lmd,\mu}$ of $\mc{G}$ with objects
	$$\ob(\mc{G}^{\lmd,\mu})=\{(A,B)\in \ob(\mc{G})|B^tA^tBA=\lmd I_n \text{ and } \tr(AB^t)=\mu\}$$
	Then $\mc{G}^{\lmd,\mu}$ is a connected cogroupoid. 
\end{lem}

Let $m,n\in \NN_{\ge 2}$ and let $A,B\in \GL_n(\kk),C,D\in \GL_m(\kk)$ such that $B^tA^tBA=\lmd I_n$, $D^tC^tDC=\lmd I_m$ and $\tr (AB^t)=\tr(CD^t)=\mu$ for some $\lmd, \mu\in \kk^{\times}$. Then $\mc{G}(A,B|C,D)$ is a nonzero algebra by the above lemma. It is a $\mc{G}(A,B)$-$\mc{G}(C,D)$-bi-Galois object by \cite[Proposition 2.8]{bi1}, the cotensor product by $\mc{G}(A,B)$-$\mc{G}(C,D)$ induces the following  equivalence of monoidal categories (Theorem \ref{thm ga} and Theorem \ref{lem 4.4}),
$$\mc{M}^{\mc{G}(A,B)}\cong^\ot \mc{M}^{\mc{G}(C,D)},\;\;\;\mc{YD}_{\mc{G}(A,B)}^{\mc{G}(A,B)}\cong^\ot \mc{YD}_{\mc{G}(C,D)}^{\mc{G}(C,D)}.$$

\subsection{The resolution}
 Let $A,B\in \GL_n(\kk)\;(n\le 2)$. Now we  construct a free Yetter-Drinfeld resolution of trivial Yetter-Drinfeld module $\kk$ over the algebra $\mc{G}(A,B)$. The foundamental $n$-dimensional $\mc{G}(A,B)$-comodule $V_{A,B}$ is defined as follows. It has a basis $v_1^{A,B}, \ldots$, $v_n^{A,B}$. The comodule action $\rho:V_{A,B}\ra V_{A,B}\ot \mc{G}(A,B)$ is defined by $\rho(v_i^{A,B})=\sum^n_{k=1}v_k^{A,B}\ot u_{ki}$. The dual vector space $V^*_{A,B}$ is a $\mc{G}(A,B)$-comodule with the comodule action $\rho:V^*_{A,B}\ra V^*_{A,B}\ot \mc{G}(A,B)$ defined by $\rho(v_i^{A,B*})=\sum^n_{j=1}v_j^{A,B*}\ot S(u_{ij})$. When there is no confusion, the basis of $V_{A,B}$ will simply be denoted by $v_1, \ldots$, $v_n$.

\begin{thm}\label{main thm}
	Let $A,B\in \GL_n(\kk)\;(n\le 2)$ such that $B^tA^tBA=\lambda I_n$ for some $\lmd\in \kk^\times$ and let $\mc{G}=\mc{G}(A,B)$. Let both $V_{A,B}$ and  $W_{A,B}$ be the foundamental $n$-dimensional $\mc{G}(A,B)$-comodules. We denote by $v_i^*\ot v_j$ and $w_i^*\ot w_j$ the respective bases of $V_{A,B}^*\ot V_{A,B}$ and $W_{A,B}^*\ot W_{A,B}$, the following complex is a free Yetter-Drinfeld resolution of $\kk$ as  a Yetter-Drinfeld moduel over $\mc{G}(A,B)$,
	$$0 \!\ra\! \kk\boxtimes\mc{G}\!\xra{\psi_4}\! (V_{A,B}^*\ot V_{A,B})\boxtimes \mc{G} \op \kk\boxtimes\mc{G}\!\xra{\psi_3} \!(V_{A,B}^*\ot V_{A,B})\boxtimes \mc{G}\op (W_{A,B}^*\ot W_{A,B})\boxtimes \mc{G}$$
	\begin{equation}\label{eq 4}\hspace{42mm}\xra{\psi_2} \kk\boxtimes\mc{G}\op(W_{A,B}^*\ot W_{A,B})\boxtimes \mc{G}\xra{\psi_1}\kk\boxtimes\mc{G}\xra{\varepsilon} \kk \ra 0, \end{equation}
	where the morphisms are defined as follows,
	$$\begin{array}{rcl}
		\psi_4:\kk\boxtimes\mc{G}&\longrightarrow &(V_{A,B}^*\ot V_{A,B})\boxtimes \mc{G} \op \kk\boxtimes\mc{G}\\
		x&\longmapsto &\sum_{i,j} v_i^*\ot v_j\ot ((A^tB)_{ij}-(AuB^t)_{ij})x+(\dd-1) x,
	\end{array}$$
	$$\psi_3:(V_{A,B}^*\ot V_{A,B})\boxtimes \mc{G} \op \kk\boxtimes\mc{G}\longrightarrow (V_{A,B}^*\ot V_{A,B})\boxtimes \mc{G}\op (W_{A,B}^*\ot W_{A,B})\boxtimes \mc{G}$$ is defined  by
	$$\begin{array}{rcl}
		\psi'_3:(V_{A,B}^*\ot V_{A,B})\boxtimes \mc{G} &\longrightarrow &(V_{A,B}^*\ot V_{A,B})\boxtimes \mc{G}\op (W_{A,B}^*\ot W_{A,B})\boxtimes \mc{G}\\
		v_i^*\ot v_j\ot x &\longmapsto & v_i^*\ot v_j\ot x +\sum_{k,l} v_k^*\ot v_l\ot (uB^t)_{il}(B^{-1})_{jk}x\\
		&&\hspace{-5mm}+\sum_{k,l}w_k^*\ot w_l\ot \frac{1}{\lmd}(B^tA^t)_{ik}(BA)_{lj}\dd x\!-\!w_i^*\ot w_j \ot x,\\
		\psi''_3:\kk\boxtimes \mc{G} &\longrightarrow &(V_{A,B}^*\ot V_{A,B})\boxtimes \mc{G}\op (W_{A,B}^*\ot W_{A,B})\boxtimes \mc{G}\\
		x&\longmapsto & \sum_{i,j}v_i^*\ot v_j\ot (A^tB)_{ij}x\\
		&&+\sum_{i,j}w_i^*\ot w_j\ot ((A uB^t)_{ij}-(A^tB)_{ij})x,
	\end{array}$$
	$$\psi_2:(V_{A,B}^*\ot V_{A,B})\boxtimes \mc{G} \op (W_{A,B}^*\ot W_{A,B})\boxtimes\mc{G}\longrightarrow \kk\boxtimes \mc{G}\op (W_{A,B}^*\ot W_{A,B})\boxtimes \mc{G}$$ is defined by
	$$\begin{array}{rcl}
		\psi'_2:(V_{A,B}^*\ot V_{A,B})\boxtimes \mc{G} &\longrightarrow &\kk\boxtimes \mc{G}\op (W_{A,B}^*\ot W_{A,B})\boxtimes \mc{G}\\
		v_i^*\ot v_j\ot x&\longmapsto & \!\!\!\!(u_{ij}-\dt_{ij})x\!\!+\!\!\sum_{k,l}w_k^*\ot w_l\ot \frac{1}{\lmd}(B^tA^t)_{ik}(BA)_{lj}\dd x\\
		&&-w_i^*\ot w_j \ot x ,\\
		\psi''_2:(W_{A,B}^*\ot W_{A,B})\boxtimes \mc{G} &\longrightarrow &\kk\boxtimes \mc{G}\op (W_{A,B}^*\ot W_{A,B})\boxtimes \mc{G}\\
		w_i^*\ot w_j\ot x&\longmapsto & -\dt_{ij}x-w_i^*\ot w_j\ot x \\&&-\sum_{k,l} w_k^*\ot w_l\ot (uB^t)_{il}(B^{-1})_{jk}x
		,
	\end{array}$$
	
	$$\psi_1:\kk\boxtimes \mc{G} \op (W_{A,B}^*\ot W_{A,B})\boxtimes\mc{G}\longrightarrow \kk\boxtimes \mc{G}$$
	is defined by
	$$\begin{array}{rcl}
		\psi'_1:\kk\boxtimes \mc{G} &\longrightarrow &\kk\boxtimes \mc{G}\\
		x&\longmapsto & (\dd-1)x ,\\
		\psi''_1:(W_{A,B}^*\ot W_{A,B})\boxtimes \mc{G} &\longrightarrow &\kk\boxtimes \mc{G}\\
		w_i^*\ot w_j\ot x&\longmapsto & (\dt_{ij}-u_{ij})x.
	\end{array}$$
\end{thm}

\subsection{Proof of Theorem \ref{main thm}} In this subsection, we keep the notations as in Theorem \ref{main thm}. 

 Let $H$ be a Hopf algebra. Let  $V\in \mc{M}^H$ and $X\in \mc{YD}^H_H$.  It is direct to verify that the following is a natural isomorphism (cf. \cite[Proposition 3.3]{bi})
 \begin{equation}\label{coli} \begin{array}{rcl} \Hom_{\mc{M}^H}(V,X)&\longrightarrow& \Hom_{\mc{YD}^H_H}(V\boxtimes H,X)\\
 		f&\longmapsto&\tilde{f}\left(v\ot h\mapsto f(v)h\right). 	
 \end{array}\end{equation}

\begin{lem}
The maps in sequence (\ref{eq 4}) in Theorem \ref{main thm}  are morphisms of Yetter-Drinfeld modules. 
\end{lem}
\proof  It is routine to check that the following maps are morphisms of $\mc{G}$-comodules. 
$$\begin{array}{rcl} \gm^V_{_1} :{V_{A,B}}^*\ot {V_{A,B}}&\longrightarrow& \kk\boxtimes\mc{G}\\
	v_i^*\ot v_j&\longmapsto&\dt_{ij},
\end{array}$$
$$\begin{array}{rcl} \gm^V_{_2} :{V_{A,B}}^*\ot {V_{A,B}}&\longrightarrow& \kk\boxtimes\mc{G}\\
	v_i^*\ot v_j&\longmapsto&u_{ij},
\end{array}$$
$$\begin{array}{rcl}\gm^V_{_3}:{V_{A,B}}^*\ot {V_{A,B}}&\longrightarrow &{V_{A,B}}^*\ot {V_{A,B}}\boxtimes \mc{G}\\
	v^*_i\ot v_j&\longmapsto&\sum_{k,l}v_k^*\ot v_l\ot (uB^t)_{il}(B^{-1})_{jk},
\end{array}$$
$$\begin{array}{rcl}\gm^V_{_4} :\kk&\longrightarrow& {V_{A,B}}^*\ot {V_{A,B}}\ot \mc{G}\\
	1&\longmapsto&\sum_{i,j} v_i^*\ot v_j\ot (A^tB)_{ij},	
\end{array}$$
$$\begin{array}{rcl} \gm^V_{_5} :\kk&\longrightarrow& {V_{A,B}}^*\ot {V_{A,B}}\ot \mc{G}\\
	1&\longmapsto&\sum_{i,j} v_i^*\ot v_j\ot (AuB^t)_{ij}.	
\end{array}$$
By replacing the space $V_{A,B}$ by the space $W_{A,B}$ in the above morphisms ,  we obtain the $\mc{G}$-comodule morphisms $\gm^W_{_i}$, $i=1,2,\cdots,5$.  Moreover, we also have  the following two morphisms of $\mc{G}$-comodules,
$$\begin{array}{rcl} \gm_{_6} :\kk&\longrightarrow&\kk\boxtimes\mc{G}\\
	1&\longmapsto&\dd-1,	
\end{array}$$
$$\begin{array}{rcl}\gm_{_7}:V_{A,B}^*\ot V_{A,B}&\longrightarrow &W_{A,B}^*\ot W_{A,B}\boxtimes \mc{G}\\
	v^*_i\ot v_j&\longmapsto&\sum_{k,l}w_k^*\ot w_l\ot \frac{1}{\lmd}(B^tA^t)_{ik}(BA)_{lj}\dd-w^*_i\ot w_j.
\end{array}$$

In sequence (\ref{eq 4}), we have $\psi'_1=\tilde{\gm_{_6}}$, $\psi''_1=\tilde{\gm^W_{_1}}-\tilde{\gm^W_{_2}}$,  $\psi'_2=\tilde{\gm^V_{_2}}-\tilde{\gm^V_{_1}}+\tilde{\gm_{_7}}$,  $\psi''_2=-\id-\tilde{\gm^W_{_1}}-\tilde{\gm^W_{_3}}$, $\psi'_3=\id+\tilde{\gm^V_{_3}}+\tilde{\gm_7}$, $\psi''_3=\tilde{\gm^W_{_5}}-\tilde{\gm^W_{_4}}+\tilde{\gm^V_{_4}}$, 
$\psi_4=\tilde{\gm^V_{_4}}-\tilde{\gm^V_{_5}}+\tilde{\gm_{_6}}$,  where the notation $\tilde{}\;$  has the same meaning as in the equation (\ref{coli}). Therefore, the morphsms $\psi_i$, $i=1,2,3,4$, are all morphisms of Yetter-Drinfeld modules. \qed

\begin{lem}
	The sequence (\ref{eq 4}) in Theorem \ref{main thm} is a complex. 
\end{lem}
\proof It is straightforward to check that $\vps\psi_1=0$. It is also straightforward to check that the following equations hold.
$$\tilde{\gm^U_{_1}}\tilde{\gm^U_{_3}}=\tilde{\gm^U_{_2}},\,\tilde{\gm^U_{_3}}\tilde{\gm^U_{_4}}=\tilde{\gm^U_{_5}},\,\tilde{\gm^U_{_3}}\tilde{\gm^U_{_5}}=\tilde{\gm^U_{_4}}+\tilde{\gm^U_{_4}}\tilde{\gm^U_{_6}},\,\tilde{\gm^U_{_2}}\tilde{\gm^U_{_3}}=\tilde{\gm_{_1}^U}+\tilde{\gm_{_6}}\tilde{\gm_{_1}^U},\,\text{ for } U=V,W$$
$$\tilde{\gm^V_{_1}}\tilde{\gm^V_{_4}}=\tilde{\gm^W_{_1}}\tilde{\gm^W_{_4}},\;\;\tilde{\gm^V_{_2}}\tilde{\gm^V_{_4}}=\tilde{\gm^W_{_1}}\tilde{\gm^W_{_5}}, \;\;\tilde{\gm_{_7}}\tilde{\gm^V_{_i}}=\tilde{\gm^W_{_i}}\tilde{\gm_{_6}},\;i=4,5,$$
$$\tilde{\gm^V_{_2}}\tilde{\gm^V_{_3}}=\tilde{\gm_{_1}^V}+\tilde{\gm_{_1}^W}\tilde{\gm_{_7}},\, \tilde{\gm_{_7}}\tilde{\gm_{_3}^V}=\tilde{\gm_{_3}^W}\tilde{\gm_{_7}},\,\tilde{\gm_{_6}}\tilde{\gm_{_2}^V}=\tilde{\gm_{_2}^W}\tilde{\gm_{_7}}.$$ 
From these equations, we obtain that $\psi_1\psi_2=\psi_2\psi_3=\psi_3\psi_4=0$. Therefore, the sequence (\ref{eq 4}) is a complex. \qed

We need some preparations to prove that the sequence (\ref{eq 4}) in  Theorem \ref{main thm} is exact.

\begin{lem}\label{lemma 2}
Let $m,n\in \NN , n,m\le 2$ and let $A,B\in \GL_n(\kk),C,D\in \GL_m(\kk)$ such that $B^tA^tBA=\lmd I_n$, $D^tC^tDC=\lmd I_m$ and $\tr (AB^t)=\tr(CD^t)=\mu$ for some $\lmd, \mu\in \kk^{\times}$. Then the sequence in Theorem \ref{main thm} is exact for $\mc{G}(A,B)$ if and only if it is exact for $\mc{G}(C,D)$. 
\end{lem}

\proof  By \cite[Theorem 4.1]{bi}, Theorem \ref{lem 4.4} and Lemma \ref{lem G E}, we obtain that the functor $$-\sq_{\mc{G}(A,B)}\mc{G}(A,B|C,D):\mc{YD}^{\mc{G(A,B)}}_{\mc{G(A,B)}}\ra \mc{YD}^{\mc{G(C,D)}}_{\mc{G(C,D)}}$$ is an equivalence of monoidal categories that preserves freeness of Yetter-Drinfeld modules. We need to check that this functor transforms the complex of Yetter-Drinfeld modules of Theorem \ref{main thm} for $\mc{G}(A,B)$ into the complex  of Yetter-Drinfeld modules of Theorem \ref{main thm} for $\mc{G}(C,D)$. 


For the simplicity, in the following, we write $V_{A,B}^*\ot V_{A,B}=\mc{V}_{A,B}$ and $W_{A,B}^*\ot W_{A,B}=\mc{W}_{A,B}$.

First, consider the following diagram. 
{\tiny$$\xymatrix{\kk \boxtimes \mc{G}(C,D)\ar[d] \ar[rr]^{\psi_4^{C,D}}&&\mc{V}_{C,D}\boxtimes \mc{G}(C,D) \op \kk\boxtimes\mc{G}(C,D)\ar[d] \\
	(\kk\sq_{\mc{G}(A,B)}\mc{G}(A,B|C,D))\boxtimes \mc{G}(C,D)\ar[d]&&((\mc{V}_{A,B}\op\kk)\sq _{\mc{G}(A,B)}\mc{G}(A,B|C,D))\boxtimes \mc{G}(C,D) \ar[d]\\	
	(\kk\boxtimes \mc{G}(A,B))\sq_{\mc{G}(A,B)}\mc{G}(A,B|C,D)\ar[rr]^{\psi_4^{A,B}\ot \id}&&(\mc{V}_{A,B}\op\kk)\boxtimes \mc{G}(A,B)\sq_{\mc{G}(A,B)}\mc{G}(A,B|C,D)
}$$}
Since $\mc{G}(A,B|C,D)$ is a left $\mc{G}(A,B)$-Galois object, we have the isomorphism $$\kk\cong \kk\sq_{\mc{G}(A,B)}\mc{G}(A,B|C,D).$$
We also have the following isomorphism of $\mc{G}(C,D)$-comodules.
$$\begin{array}{rcl}
	\mc{V}_{C,D}&\lra& 	\mc{V}_{A,B}\sq_{\mc{G}(A,B)}\mc{G}(A,B|C,D),\\
	v_i^{C,D*}\ot v_j^{C,D}&\longmapsto&\sum_{k,l} v_k^{A,B*}\ot v_l^{A,B}\ot S_{CD,AB}(u_{ik}^{CD,AB})u_{lj}^{AB,CD}.
\end{array}$$
Together with the isomorphism of Theorem \ref{lem 4.4}, we obtain that  the composition 
{\small $$\kk \boxtimes \mc{G}(C,D)\ra(\kk\sq_{\mc{G}(A,B)}\mc{G}(A,B|C,D))\boxtimes \mc{G}(C,D)\ra(\kk\boxtimes \mc{G}(A,B))\sq_{\mc{G}(A,B)}\mc{G}(A,B|C,D)$$}
is given by
$$\begin{array}{rcl}\kk \boxtimes \mc{G}(C,D)&\lra&(\kk\boxtimes \mc{G}(A,B))\sq_{\mc{G}(A,B)}\mc{G}(A,B|C,D)\\
	1\ot x^{CD,CD}&\longmapsto& x_{(2)}^{AB,AB}\ot S_{CD,AB}(x_{(1)}^{CD,AB})x_{(3)}^{AB,CD},
\end{array}$$
and the composition 
{\tiny$$\mc{V}_{C,D}\boxtimes \mc{G}(C,D) \ra (\mc{V}_{A,B}\sq _{\mc{G}(A,B)}\mc{G}(A,B|C,D))\boxtimes \mc{G}(C,D)\ra \mc{V}_{A,B}\boxtimes \mc{G}(A,B)\sq_{\mc{G}(A,B)}\mc{G}(A,B|C,D)  $$}
is given by 
{\tiny $$\begin{array}{rcl}\mc{V}_{C,D}\boxtimes \mc{G}(C,D)&\ra&\mc{V}_{A,B}\boxtimes \mc{G}(A,B)\sq_{\mc{G}(A,B)}\mc{G}(A,B|C,D) \\
	{\tiny v_i^{C,D*}\ot v_j^{C,D}\ot x^{CD,CD}}&\mapsto& \sum_{k,l} v_k^{A,B*}\ot v_l^{A,B}\\
	&&\ot x_{(2)}^{AB,AB}\ot S_{CD,AB}(x_{(1)}^{CD,AB}) S_{CD,AB}(u_{ik}^{CD,AB})u^{AB,CD}_{lj}x_{(3)}^{AB,CD}.
\end{array}$$}
The vertical arrows are compositions of isomorphisms, so are isomorphisms.  It is routine but tedious to check that  the above diagram is commutative. Similarly, we have the following commutative diagrams with vertical morphisms are isomorphisms. 

{\tiny$$\xymatrix{\mc{V}_{C,D}\boxtimes \mc{G}(C,D) \op \kk\boxtimes\mc{G}(C,D)\ar[d] \ar[r]^{\psi_3^{C,D}}&\mc{V}_{C,D}\boxtimes \mc{G}(C,D) \op \mc{W}_{C,D}\boxtimes \mc{G}(C,D)\ar[d] \\
		((\mc{V}_{A,B}\op\kk)\sq_{\mc{G}(A,B)}\mc{G}(A,B|C,D))\boxtimes \mc{G}(C,D)\ar[d]&((\mc{V}_{A,B}\op\mc{W}_{A,B})\sq _{\mc{G}(A,B)}\mc{G}(A,B|C,D))\boxtimes \mc{G}(C,D) \ar[d]\\	
		((\mc{V}_{A,B}\op\kk)\boxtimes \mc{G}(A,B))\sq_{\mc{G}(A,B)}\mc{G}(A,B|C,D)\ar[r]^{\psi_3^{A,B}\ot \id}&((\mc{V}_{A,B}\op\mc{W}_{A,B}))\boxtimes \mc{G}(A,B)\sq_{\mc{G}(A,B)}\mc{G}(A,B|C,D),
	}$$}
{\tiny$$\xymatrix{\mc{V}_{C,D}\boxtimes \mc{G}(C,D) \op \mc{W}_{C,D}\boxtimes \mc{G}(C,D)\ar[d] \ar[r]^{\psi_2^{C,D}}& \mc{W}_{C,D}\boxtimes \mc{G}(C,D)\op \kk\boxtimes \mc{G}(C,D)\ar[d] \\
		((\mc{V}_{A,B}\op\mc{W}_{A,B})\sq_{\mc{G}(A,B)}\mc{G}(A,B|C,D))\boxtimes \mc{G}(C,D)\ar[d]&((\mc{W}_{A,B}\op \kk)\sq _{\mc{G}(A,B)}\mc{G}(A,B|C,D))\boxtimes \mc{G}(C,D) \ar[d]\\	
		((\mc{V}_{A,B}\op\mc{W}_{A,B})\boxtimes \mc{G}(A,B))\sq_{\mc{G}(A,B)}\mc{G}(A,B|C,D)\ar[r]^{\psi_2^{A,B}\ot \id}&((\mc{W}_{A,B})\op\kk)\boxtimes \mc{G}(A,B)\sq_{\mc{G}(A,B)}\mc{G}(A,B|C,D),
	}$$}
{\tiny$$\xymatrix{\mc{W}_{C,D}\boxtimes \mc{G}(C,D)\op \kk\boxtimes \mc{G}(C,D)\ar[d] \ar[r]^{\hspace{10mm}\psi_1^{C,D}}& \kk\boxtimes\mc{G}(C,D) \ar[d] \\
		((\mc{W}_{A,B}\op \kk)\sq_{\mc{G}(A,B)}\mc{G}(A,B|C,D))\boxtimes \mc{G}(C,D)\ar[d]&(\kk\sq _{\mc{G}(A,B)}\mc{G}(A,B|C,D))\boxtimes \mc{G}(C,D) \ar[d]\\	
		((\mc{W}_{A,B}\op \kk)\boxtimes \mc{G}(A,B))\sq_{\mc{G}(A,B)}\mc{G}(A,B|C,D)\ar[r]^{\hspace{8mm}\psi_1^{A,B}\ot\id}&\kk\boxtimes \mc{G}(A,B)\sq_{\mc{G}(A,B)}\mc{G}(A,B|C,D).
	}$$}Thus we obtain that the sequence in Theorem \ref{main thm} is exact for $\mc{G}(A,B)$ if and only if it is exact for $\mc{G}(C,D)$. \qed

In order to show that the sequence (\ref{eq 4}) in Theorem \ref{main thm} is exact, we first show that the sequence is exact when $A=A_q, B=A_q^{-1}$. It takes several steps. We know that $\mc{G}(A_q,A_q^{-1})=\mc{O}(\GL_q(2))$. As usual, we put $a=u_{11}$, $b=u_{12}$, $c=u_{21}$ and $d=u_{22}$. The element $\dd$ is the just the quantum determinant $\dd_q:=ad-qbc=da-q^{-1}bc$, which is an central element.

Recall that $\mc{O}(\SL_q(2))$ for some $q\in \kk^\times$ is the algebra generated by $\bar{a},\bar{b},\bar{c},\bar{d}$, subject to the relations 
$$\begin{array}{ccc}
	\bar{a}\bar{b}=q\bar{b}\bar{a}\;\;&\;\;\bar{a}\bar{c}=q\bar{c}\bar{a}\;\;&\;\;\bar{b}\bar{c}=\bar{c}\bar{b}\\
	\bar{b}\bar{d}=q\bar{d}\bar{b}\;\;&\;\;\bar{c}\bar{d}=q\bar{d}\bar{c}\;\;&\;\;\bar{a}\bar{d}-q\bar{b}\bar{c}=\bar{d}\bar{a}-q^{-1}\bar{b}\bar{c}=1
\end{array}$$

\begin{lem}\label{lem: iso}There is an algebra automorphism
	$$\mc{O}(\GL_q(2))\cong 	\mc{O}(\SL_q(2))[z^{\pm 1}].$$
\end{lem}
\proof It is routine to check that the map from $\mc{O}(\GL_q(2))$ to $\mc{O}(\SL_q(2))[z^{\pm 1}]$ mapping $a,b,c,d$ to $\bar{a}z,\bar{b}z,\bar{c},\bar{d}$, respectively, is an algebra automorphism. 

We point out regarding the notions in \cite{bi}, we have $E_q=A_{q^{-1}}$. So $\mathcal B(E_q)=\mathcal O(SL_{q^{-1}}(2))$ as in our paper. Then following result is Lemma 5.6 in \cite{bi} in terms of $O(SL_{q}(2))=B(E_{q^{-1}})$ in light of \cite[Lemma 5.5]{bi}.

\begin{lem}\label{lem1}
	Let $\mc{A}=\mc{O}(\SL_q(2))$, $V$ be a 2-dimensional vector space with a basis $\{v_1,v_2\}$.  There exists a free resolution of $\kk$ as an $\mc{A}$-module, 
\begin{equation}\label{complex} 0 \ra \mc{A}\xra{\phi_3} (V^*\ot V)\ot \mc{A} \xra{\phi_2} (V^*\ot V)\ot \mc{A} \xra{\phi_1}\mc{A}\xra{\vps} \kk\ra 0.\end{equation}
	The morphisms $\phi_1, \phi_2,\phi_3$ are defined as follows,

$\phi_3(x)=v_1^*\ot v_1\ot ((-q+q^{-1}d)x)+v_1^*\ot v_2\ot (-cx)+v_2^*\ot v_1\ot (-bx)+v_2^*\ot v_2\ot ((-q^{-1}+qa)x),$

$\phi_2(v_1^*\ot v_1\ot x)=v_1^*\ot v_1\ot x +v_2^*\ot v_1\ot (-q^{-1}bx) +v_2^*\ot v_2\ot ax, $

$\phi_2(v_1^*\ot v_2\ot x)=v_1^*\ot v_1\ot bx +v_1^*\ot v_2\ot (1-qa)x,$

$\phi_2(v_2^*\ot v_1\ot x)=v_2^*\ot v_1\ot (1-q^{-1}d)x +v_2^*\ot v_2\ot cx,$

$\phi_2(v_2^*\ot v_2\ot x)=v_1^*\ot v_1\ot dx +v_1^*\ot v_2\ot (-qcx) +v_2^*\ot v_2\ot x,$

$\phi_1(v_1^*\ot v_1\ot x)=(a-1)x$, $\phi_1(v_1^*\ot v_2\ot x)=bx,$

$\phi_1(v_2^*\ot v_1\ot x)=cx$, $\phi_1(v_2^*\ot v_2\ot x)=(d-1)x.$
\end{lem}

Applying the functor $-\ot_{\mc{A}}\mc{A}[z^\pm 1]$ to the complex (\ref{complex}), we obtain the following complex,
$$0\ra \mc{A}[z^{\pm 1}]\xra{  \phi'_3} (V^*\ot V)\ot \mc{A}[z^{\pm 1}]\xra{  \phi'_2}(V^*\ot V)\ot \mc{A}[z^{\pm 1}]\xra{\phi'_1}  \mc{A}[z^{\pm 1}],$$
where $\phi'_i=\phi_i\ot_{\mc{A}} \mc{A}[z^{\pm 1}]$ $ (i=1,2,3)$.
Then we have the following commutative diagram, 
\begin{equation}{\tiny\label{eq 1}
		\xymatrix{0\ar[r] & \mc{A}[z^{\pm 1}]\ar[d]^{f_{_3}}\ar[r]^{  \phi'_3\;\;\;\;\;\;\;\;\;\;\;\;}& (V^*\ot V)\ot \mc{A}[z^{\pm 1}]\ar[d]^{f_{_2}}\ar[r]^{\phi'_2}& (V^*\ot V)\ot \mc{A}[z^{\pm 1}]\ar[d]^{f_{_1}}\ar[r]^{\;\;\;\;\;\;\;\;\;\;\;\;\phi'_1}& \mc{A}[z^{\pm 1}]\ar[d]^{f_{_0}}\\
			0\ar[r] & \mc{A}[z^{\pm 1}]\ar[r]^{\phi'_3\;\;\;\;\;\;\;\;\;\;\;\;}& (V^*\ot V)\ot \mc{A}[z^{\pm 1}]\ar[r]^{\phi'_2}& (V^*\ot V)\ot \mc{A}[z^{\pm 1}]\ar[r]^{\;\;\;\;\;\;\;\;\;\;\;\;\phi'_1}& \mc{A}[z^{\pm 1}]
	}}
\end{equation}
where  $f_i \;(i=0,1,2,3)$ is the morphism of the right multiplication by $z-1$.
\begin{lem}\label{res}
	We have 
$$cone(f)\ra  \kk\ra 0$$
is a projective resolution of $\kk$ as an $\mc{A}[z^{\pm 1}]$-module. 
\end{lem}

\proof  
 Let $C_1$ be the following complex 
$$0\ra \kk  \ot_{\mc{A}}  \mc{A}[z^{\pm 1}]\xra{\cdot (z-1)}\kk  \ot_{\mc{A}}  \mc{A}[z^{\pm 1}]\ra 0.$$
We have a morphism $g$ from $cone(f)$ to this complex,
{\tiny$$\xymatrix{0\ar[r] & \mc{A}[z^{\pm 1}]\ar[d]^{-f_3}\ar[r]^{  -\phi'_3\;\;\;\;\;\;\;\;\;\;\;\;}& (V^*\!\!\ot\!\! V)\!\ot\! \mc{A}[z^{\pm 1}]\ar[d]^{-f_2}\ar[r]^{-\phi'_2}& (V^*\!\!\ot\!\! V)\!\ot\! \mc{A}[z^{\pm 1}]\ar[d]^{-f_1}\ar[r]^{\;\;\;\;\;\;\;\;\;\;\;\;-\phi'_1}& \mc{A}[z^{\pm 1}]\ar[d]^{-f_0}\ar[r]^{g_{_1} }&\kk \!\!\ot_{\mc{A}} \!\!\mc{A}[z^{\pm 1}]\ar[d]^{\cdot (z-1)}\\
		0\ar[r] & \mc{A}[z^{\pm 1}]\ar[r]^{\phi'_3\;\;\;\;\;\;\;\;\;\;\;\;}& (V^*\!\!\ot\!\! V)\!\ot\! \mc{A}[z^{\pm 1}]\ar[r]^{\phi'_2}& (V^*\!\!\ot\!\! V)\!\ot\! \mc{A}[z^{\pm 1}]\ar[r]^{\;\;\;\;\;\;\;\;\;\;\;\;\phi'_1}& \mc{A}[z^{\pm 1}]\ar[r]^{g_{_0}}& \kk \!\!\ot_{\mc{A}} \!\!\mc{A}[z^{\pm 1}],
	}$$}
where $g_{_0}=\varepsilon\ot \mc{A}[z^{\pm 1}]$ and $g_{_1}=-\varepsilon\ot \mc{A}[z^{\pm 1}]$.

Since the functor $-\ot_{\mc{A}} \mc{A}[z^{\pm 1}]$ is exact, it is easy to see that we have the following bicomplex with exact rows.
{\tiny$$\xymatrix{0\ar[r] & \mc{A}[z^{\pm 1}]\ar[d]^{f_3}\ar[r]^{  \phi'_3\;\;\;\;\;\;\;\;\;\;\;\;}& (V^*\!\!\ot\!\! V)\!\ot\! \mc{A}[z^{\pm 1}]\ar[d]^{f_2}\ar[r]^{\phi'_2}& (V^*\!\!\ot\!\! V)\!\ot\! \mc{A}[z^{\pm 1}]\ar[d]^{f_1}\ar[r]^{\;\;\;\;\;\;\;\;\;\;\;\;\phi'_1}& \mc{A}[z^{\pm 1}]\ar[d]^{f_0}\ar[r]^{\varepsilon\ot \mc{A}[z^{\pm 1}] }&\kk \!\!\ot_{\mc{A}} \!\!\mc{A}[z^{\pm 1}]\ar[d]^{\cdot (z-1)}\\
	0\ar[r] & \mc{A}[z^{\pm 1}]\ar[r]^{-\phi'_3\;\;\;\;\;\;\;\;\;\;\;\;}& (V^*\!\!\ot\!\! V)\!\ot\! \mc{A}[z^{\pm 1}]\ar[r]^{-\phi'_2}& (V^*\!\!\ot\!\! V)\!\ot\! \mc{A}[z^{\pm 1}]\ar[r]^{\;\;\;\;\;\;\;\;\;\;\;\;-\phi'_1}& \mc{A}[z^{\pm 1}]\ar[r]^{-\varepsilon\ot \mc{A}[z^{\pm 1}] }& \kk \!\!\ot_{\mc{A}} \!\!\mc{A}[z^{\pm 1}]. 
}$$}Denote the above bicomplex by $C_2$. By the Acyclic Assembly Lemma \cite[Lemma 2.7.3]{we},  the total complex ${\rm Tot}(C_2)$ is exact.  The total complex ${\rm Tot}(C_2)$ is just the complex $cone(g)$.  Thus, we have  $cone(f)$ and $C_1$ are quasi-isomorphic. We have the following exact sequence of $\mc{A}[z^{\pm 1}]$-modules. 
$$0\ra \kk \ot_{\mc{A}} \mc{A}[z^{\pm 1}]\xra{\cdot(z-1)}\kk \ot_{\mc{A}} \mc{A}[z^{\pm 1}]\ra \kk\ra 0$$
Therefore, the complex $cone(f)\ra  \kk\ra 0$ is exact. The functor $-\ot_{\mc{A}} \mc{A}[z^{\pm 1}]$ also preserve projectivity, so $cone(f)\ra  \kk\ra 0$ is a projective resolution of $\kk$. \qed


\begin{lem}\label{lemma 1}
	Let $q\in \kk^\times$.  The sequence in Theorem \ref{main thm}  is exact when $A=A_q$ and $B=A_q^{-1}$.  
\end{lem}
\proof 	Let $\mc{B}=\mc{G}(A_q,A_q^{-1})$.   When $A=A_q$ and $B=A_q^{-1}$, the sequence (\ref{eq 4}) is just the following complex,  
	\begin{equation}\label{eq 2}
	0 \ra P_4\xra{\psi_4} P_3\xra{\psi_3} P_2\xra{\psi_2} P_1\xra{\psi_1} P_0 \xra{\varepsilon}\kk\ra 0,
	\end{equation}
	where $P_4=\kk\boxtimes\mc{B}$, $P_3=(V_{A,B}^*\ot V_{A,B})\boxtimes \mc{B}\op \kk\boxtimes\mc{B}$, $P_2=(V_{A,B}^*\ot V_{A,B})\boxtimes \mc{B}\op (W_{A,B}^*\ot W_{A,B})\boxtimes \mc{B}$, $P_1=(W_{A,B}^*\ot W_{A,B})\boxtimes \mc{B}\op \kk\boxtimes\mc{B}$, $P_0=\kk\boxtimes\mc{B}$ and the morphisms are given by
	
	$\psi_4(x)=v_1^*\ot v_1\ot((-q+q^{-1}d)x)+v^*_1\ot v_2\ot (-cx)+v_2^*\ot v_1\ot (-bx)+v_2^*\ot v_2\ot ((-q^{-1}+qa)x)+(\dd_q-1) x$,
	\vspace{3mm}
	
	$\psi_3(v_1^*\ot v_1\ot x)=v_1^*\ot v_1\ot x+v_2^*\ot v_1\ot(-q^{-1}bx)+v_2^*\ot v_2\ot ax+w_1^*\ot w_1\ot (\dd_q-1) x$,
	
	$\psi_3(v_1^*\ot v_2\ot x)=v_1^*\ot v_1\ot bx+v_1^*\ot v_2\ot (1-qa)x+w_1^*\ot w_2\ot (\dd_q-1) x$,
	
	$\psi_3(v_2^*\ot v_1\ot x)=v_2^*\ot v_1\ot (1-q^{-1}d)x+v_2^*\ot v_2\ot cx+w_2^*\ot w_1\ot (\dd_q-1) x$,
	
	$\psi_3(v_2^*\ot v_2\ot x)=v_1^*\ot v_1\ot dx+v_1^*\ot v_2\ot(-qc)x+v_2^*\ot v_2\ot x+w_2^*\ot w_2\ot (\dd_q-1) x$,
	
	$\psi_3(x)=w_1^*\ot w_1\ot (q-q^{-1}d)x+w_1^*\ot w_2\ot cx +w_2^*\ot w_1\ot bx+w_2^*\ot w_2\ot (q^{-1}-qa)x-v_1^*\ot v_1\ot qx-v_2^*\ot v_2\ot q^{-1}x$,
	\vspace{3mm}
	
	$\psi_2(v_1^*\ot v_1\ot x)=(a-1)x+w_1^*\ot w_1\ot (\dd_q-1) x$,
	
	$\psi_2(v_1^*\ot v_2\ot x)=bx+w_1^*\ot w_2\ot (\dd_q-1) x$,
	
	$\psi_2(v_2^*\ot v_1\ot x)=cx+w_2^*\ot w_1\ot (\dd_q-1) x$,
	
	$\psi_2(v_2^*\ot v_2\ot x)=(d-1)x+w_2^*\ot w_2\ot (\dd_q-1) x$,
	
	$\psi_2(w_1^*\ot w_1\ot x)=-w^*_1\ot w_1\ot x-w^*_2\ot w_1\ot (-q^{-1}b)x-w^*_2\ot w_2\ot ax-x$,
	
	$\psi_2(w_1^*\ot w_2\ot x)=-w^*_1\ot w_1\ot bx-w^*_1\ot w_2\ot (1-qa)x$,
	
	$\psi_2(w_2^*\ot w_1\ot x)=-w^*_2\ot w_1\ot (1-q^{-1}d)x-w^*_2\ot w_2\ot cx$,
	
	$\psi_2(w_2^*\ot w_2\ot x)=-w^*_1\ot w_1\ot dx-w^*_1\ot w_2\ot (-qc)x-w^*_2\ot w_2\ot x-x$,
	\vspace{3mm}
	
	$\psi_1(w_1^*\ot w_1\ot x)=(1-a)x$, 	$\psi_1(w_1^*\ot w_2\ot x)=-bx$,
	
	$\psi_1(w_2^*\ot w_1\ot x)=-cx$, 	$\psi_1(w_2^*\ot w_2\ot x)=(1-d)x$, 	$\psi_1(x)=(\dd_q-1)x$.

Applying the algebra isomorphism $\mc{O}(\GL_q(2))\cong\mc{O}(\SL_q(2))[z^{\pm 1}]$ stated in Lemma \ref{lem: iso} and Lemma \ref{res}, we obtain that the following complex is exact,
\begin{equation}\label{eq 3}
		0 \ra P_4\xra{\bar{\psi}_4} P_3\xra{\bar{\psi}_3} P_2\xra{\bar{\psi}_2} P_1\xra{\bar{\psi}_1} P_0 \xra{\varepsilon}\kk\ra 0,
\end{equation}
where $\bar{\psi}_i$, $i=1,2,3,4$, are morphisms of $\mc{B}$-modules defined as follows,

$\bar{\psi}_4(x)=v_1^*\ot v_1\ot((-q+q^{-1}d)x)+v^*_1\ot v_2\ot (-cx)+v_2^*\ot v_1\ot (-b\dd_q^{-1}x)+v_2^*\ot v_2\ot ((-q^{-1}+qa\dd_q^{-1})x)+(\dd_q-1) x$,
\vspace{3mm}

$\bar{\psi}_3(v_1^*\ot v_1\ot x)=v_1^*\ot v_1\ot x+v_2^*\ot v_1\ot(-q^{-1}b\dd_q^{-1}x)+v_2^*\ot v_2\ot a\dd_q^{-1}x+w_1^*\ot w_1\ot (\dd_q-1) x$,

$\bar{\psi}_3(v_1^*\ot v_2\ot x)=v_1^*\ot v_1\ot b\dd_q^{-1}x+v_1^*\ot v_2\ot (1-qa\dd_q^{-1})x+w_1^*\ot w_2\ot (\dd_q-1) x$,

$\bar{\psi}_3(v_2^*\ot v_1\ot x)=v_2^*\ot v_1\ot (1-q^{-1}d)x+v_2^*\ot v_2\ot cx+w_2^*\ot w_1\ot (\dd_q-1) x$,

$\bar{\psi}_3(v_2^*\ot v_2\ot x)=v_1^*\ot v_1\ot dx+v_1^*\ot v_2\ot(-qc)x+v_2^*\ot v_2\ot x+w_2^*\ot w_2\ot (\dd_q-1) x$,

$\bar{\psi}_3(x)=w_1^*\ot w_1\ot (q-q^{-1}d)x+w_1^*\ot w_2\ot cx +w_2^*\ot w_1\ot b\dd_q^{-1}x+w_2^*\ot w_2\ot (q^{-1}-qa\dd_q^{-1})x$,
\vspace{3mm}

$\bar{\psi}_2(v_1^*\ot v_1\ot x)=(a\dd_q^{-1}-1)x+w_1^*\ot w_1\ot (\dd_q-1) x$,

$\bar{\psi}_2(v_1^*\ot v_2\ot x)=b\dd_q^{-1}x+w_1^*\ot w_2\ot (\dd_q-1) x$,

$\bar{\psi}_2(v_2^*\ot v_1\ot x)=cx+w_2^*\ot w_1\ot (\dd_q-1) x$,

$\bar{\psi}_2(v_2^*\ot v_2\ot x)=(d-1)x+w_2^*\ot w_2\ot (\dd_q-1) x$,

$\bar{\psi}_2(w_1^*\ot w_1\ot x)=-w^*_1\ot w_1\ot x-w^*_2\ot w_1\ot (-q^{-1}b\dd_q^{-1})x-w^*_2\ot w_2\ot a\dd_q^{-1}x$,

$\bar{\psi}_2(w_1^*\ot w_2\ot x)=-w^*_1\ot w_1\ot b\dd_q^{-1}x-w^*_1\ot w_2\ot (1-qa\dd_q^{-1})x$,

$\bar{\psi}_2(w_2^*\ot w_1\ot x)=-w^*_2\ot w_1\ot (1-q^{-1}d)x-w^*_2\ot w_2\ot cx$,

$\bar{\psi}_2(w_2^*\ot w_2\ot x)=-w^*_1\ot w_1\ot dx-w^*_1\ot w_2\ot (-qc)x-w^*_2\ot w_2\ot x$,
\vspace{3mm}

$\bar{\psi}_1(w_1^*\ot w_1\ot x)=(1-a\dd_q^{-1})x$, 	$\bar{\psi}_1(w_1^*\ot w_2\ot x)=-b\dd_q^{-1}x$,

$\bar{\psi}_1(w_2^*\ot w_1\ot x)=-cx$, 	$\bar{\psi}_1(w_2^*\ot w_2\ot x)=(1-d)x$, 	$\bar{\psi}_1(x)=(\dd_q-1)x$.


We will show that complexes (\ref{eq 2}) and  (\ref{eq 3}) are isomorphic as complexes over $\mc{B}$-modules, then the complex (\ref{eq 2}) is exact. 
 
We have the following commutative diagram
$$\xymatrix
{0 \ar[r]& \ar^{g_{_4}}[d]P_4\ar^-{\psi_4}[r]&\ar^{g_{_3}}[d] P_3\ar^-{\psi_3}[r]&\ar^{g_{_2}}[d] P_2\ar^-{\psi_2}[r] &\ar^{g_{_1}}[d]P_1\ar^-{\psi_1}[r]&\ar^{g_{_0}}[d]P_0\ar[r]^{\varepsilon}&\kk\ar[d]\ar[r]&0\\
	0 \ar[r]& P_4\ar^-{\bar{\psi}_4}[r]& P_3\ar^-{\bar{\psi}_3}[r]& P_2\ar^-{\bar{\psi}_2}[r] &P_1\ar^-{\bar{\psi}_1}[r]&P_0\ar[r]^{\varepsilon}&\kk\ar[r]&0},$$
where the morphisms $g_i$, $i=0,1,2,3,4$ are defined as follows, $g_{_0}$ is just the identity, $g_{_4}(x)=x\dd_q$, for $x\in\mc{B}$, $g_{_1}, g_{_2}, g_{_3}$ are the morphisms of right $\mc{B}$-modules defined by the following maps (we also denote by $g_{_1}, g_{_2}, g_{_3}$ respectively as well),
{\small $$ \begin{array}{rl}	
&g_{_1}\begin{pmatrix}
	w_1^*\ot w_1,\!\!\!&\!\!\!  w_1^*\ot w_2,\!\!\!&\!\!\!  w_2^*\ot w_1,\!\!\!&\!\!\!  w_2^*\ot w_2,\!\!\!&\!\!\! 1
	\end{pmatrix}\\
=&\begin{pmatrix}
	w_1^*\ot w_1\ot 1,\!\!\!&\!\!\!  w_1^*\ot w_2\ot 1,\!\!\!&\!\!\!  w_2^*\ot w_1\ot 1,\!\!\!&\!\!\!  w_2^*\ot w_2\ot 1,\!\!\!&\!\!\! 1
\end{pmatrix}\begin{pmatrix}
\dd_q
& & & &\\
&\dd_q & & &\\
&  & 1 & &\\
&  & & 1 &\\
-1&  & &  &1\\
\end{pmatrix},
	\end{array}
$$}
{\tiny $$ \begin{array}{rl}	
		\!\!\!&\!\!\!g_{_2}\begin{pmatrix}
				v_1^*\ot v_1,\!\!\!&\!\!\!  v_1^*\ot v_2,\!\!\!&\!\!\!  v_2^*\ot v_1,\!\!\!&\!\!\!  v_2^*\ot v_2,\!\!\!&\!\!\! w_1^*\ot w_1,\!\!\!&\!\!\! w_1^*\ot w_2,\!\!\!&\!\!\! w_2^*\ot w_1,\!\!\!&\!\!\!  w_2^*\ot w_2
		\end{pmatrix}\\
		=\!\!\!&\!\!\!\begin{pmatrix}
			v_1^*\ot v_1\ot 1,\!\!\!&\!\!\!  v_1^*\ot v_2\ot 1,\!\!\!&\!\!\!  v_2^*\ot v_1\ot 1,\!\!\!&\!\!\!  v_2^*\ot v_2\ot 1,\!\!\!&\!\!\! w_1^*\ot w_1\ot 1,\!\!\!&\!\!\! w_1^*\ot w_2\ot 1,\!\!\!&\!\!\! w_2^*\ot w_1\ot 1,\!\!\!&\!\!\!  w_2^*\ot w_2\ot 1
		\end{pmatrix}P,
	\end{array}
	$$}
where $P$ is the matirx {\tiny $$\begin{pmatrix}
	\dd_q	& & & & & & &\\
	&\dd_q & & & &\dd_q & &\\
	&  & 1 & & & & &\\
	&  & & 1 & & & & 1\\
	&  & &  &\dd_q & & &\\
	&  & &  & & \dd_q^2  & &\\
	&  & &  & &  &1 &\\
	&  & &  & &  &  &\dd_q\\
\end{pmatrix},$$}
{\small $$ \begin{array}{rl}	
		&g_{_3}\begin{pmatrix}
			v_1^*\ot v_1,\!\!\!&\!\!\!  v_1^*\ot v_2,\!\!\!&\!\!\!  v_2^*\ot v_1,\!\!\!&\!\!\!  v_2^*\ot v_2,\!\!\!&\!\!\! 1
		\end{pmatrix}\\
		=&\begin{pmatrix}
			v_1^*\ot v_1\ot 1,\!\!\!&\!\!\!  v_1^*\ot v_2\ot 1,\!\!\!&\!\!\!  v_2^*\ot v_1\ot 1,\!\!\!&\!\!\!  v_2^*\ot v_2\ot 1,\!\!\!&\!\!\! 1
		\end{pmatrix}\begin{pmatrix}
			\dd_q
			& & & &\\
			&\dd^2_q & & &c\dd_q\\
			&  & 1 & &\\
			&  & & \dd_q &-qa\\
			&  & &  &\dd_q\\
		\end{pmatrix}.
	\end{array}
	$$}It is routine to check that $g_{_i}$, $i=0,1,2,3,4$ are isomorphisms. So we conclude that the complex (\ref{eq 2}) is exact.  \qed

Now we can prove that the sequence (\ref{eq 4}) in Theorem \ref{main thm} is exact to complete the proof of Theorem \ref{main thm}. Let $A,B\in \GL_n(\kk)\;(n\le 2)$ such that $B^tA^tBA=\lambda I_n$ for some $\lmd\in \kk^\times$. Let $q\in\kk^\times$ such that $q^2+\sqrt{\lmd^{-1}}\tr(AB^t)q+1=0$. Note that $\mc{G}(A,B)=\mc{G}(a A,b B)$ for all $a, b \in \kk^\times$. Put $A'=\sqrt{\lmd^{-1}}A$ and $B'=B$. Then $B'^tA'^tB'A'= I_n$ and $\tr(A'B'^t)=\frac{1}{\sqrt{\lmd^{-1}}}\tr(AB^t)=-q^{-1}-q=\tr(A_q(A_q^{-1})^t)$. Therefore, by Lemma \ref{lemma 2} and Lemma \ref{lemma 1} and $\mc{G}(A,B)=\mc{G}(A', B')$, the sequence (\ref{eq 4}) in Theorem \ref{main thm} is exact.  
\qed


\section{The Calabi-Yau property}\label{3}
In this section, we show that the algebras $\mc{G}(A,B)$ are twisted CY algebras. 
\begin{thm}\label{thm main1}
Let $A,B\in \GL_n(\kk)\;(n\le 2)$ such that $B^tA^tBA=\lambda I_n$ for some $\lmd\in \kk^\times$. The algebra $\mc{G}=\mc{G}(A,B)$ is a twisted CY algebra with Nakayama automorphism $\mu$ defined by $\mu(u)= (A^t)^{-1}A u B^tB^{-1}$ and $\mu(\dd^{\pm 1})=\dd^{\pm 1}$.
\end{thm}

\proof We keep the notations as in Theorem \ref{main thm}. From Theorem \ref{main thm}, we see that the right trivial module $\kk$ has a bounded projective resolution with each term finitely generated. By Lemma \ref{lem AS}, in order to show that $\mc{G}(A,B)$ is twisted CY, we only need to prove that  $\mc{G}(A,B)$ is right AS-Gorenstein. 

After applying the functor $\Hom_{\mc{G}}(-,\mc{G})$ to the resolution (\ref{eq 4}) and standard identifications, we get the following complex
$$0 \ra \mc{G}\xra{\psi^t_1}  \mc{G}\ot(W_{A,B}^*\ot W_{A,B}) \op \mc{G}\xra{\psi^t_2} \mc{G}\ot(W_{A,B}^*\ot W_{A,B})\op \mc{G}\ot (V_{A,B}^*\ot V_{A,B})\hspace{10mm}$$
$$\hspace{60mm}\xra{\psi^t_3} \mc{G}\op\mc{G}\ot (V_{A,B}^*\ot V_{A,B})\xra{\psi^t_4}\mc{G}\ra 0,$$
with 
$$\begin{array}{rcl}
	\psi^t_1:\mc{G}&\longrightarrow &\mc{G}\ot(W_{A,B}^*\ot W_{A,B})  \op \mc{G}\\
	x&\longmapsto &\sum_{i,j}  x(\dt_{ji}-u_{ji})\ot w_i^*\ot w_j+x(\dd-1),
\end{array}$$
$$\psi_2^t:  \mc{G}\ot(W_{A,B}^*\ot W_{A,B}) \op \mc{G} \longrightarrow  \mc{G}\ot(W_{A,B}^*\ot W_{A,B})\op \mc{G}\ot(V_{A,B}^*\ot V_{A,B})$$ is defined  by
$$\begin{array}{rcl}
	\mc{G}\ot(W_{A,B}^*\ot W_{A,B})  &\longrightarrow & \mc{G}\ot(W_{A,B}^*\ot W_{A,B})\op \mc{G}\ot(V_{A,B}^*\ot V_{A,B})\\
	x\ot w_i^*\ot w_j &\longmapsto &  -x\ot w_i^*\ot w_j-x\sum_{k,l} (uB^t)_{li}(B^{-1})_{kj} \ot w_k^*\ot w_l\\
	&&\hspace{2mm}+x\sum_{k,l}\frac{1}{\lmd}(BA)_{ik}(B^tA^t)_{lj}\dd\ot v_k^*\ot v_l-x\ot v^*_i\ot v_j,\\
	
	\mc{G} &\longrightarrow &\mc{G}\ot(W_{A,B}^*\ot W_{A,B})\op \mc{G}\ot(V_{A,B}^*\ot V_{A,B})\\
	x&\longmapsto & -x\dt_{ij}\ot w^*_i\ot w_j+ x(u_{ji}-\dt_{ji})\ot v^*_i\ot v_j,\\
\end{array}$$
$$\psi^t_3: \mc{G}\ot(W_{A,B}^*\ot W_{A,B}) \op \mc{G}\ot(V_{A,B}^*\ot V_{A,B})\longrightarrow \mc{G}\op \mc{G}\ot(V_{A,B}^*\ot V_{A,B})$$ is defined by
$$\begin{array}{rcl}
	\mc{G}\ot(W_{A,B}^*\ot W_{A,B}) &\longrightarrow &\mc{G}\op \mc{G}\ot(V_{A,B}^*\ot V_{A,B})\\
	x\ot w_i^*\ot w_j&\longmapsto & x((AuB^t)_{ji}-(A^tB)_{ji})\\&&+x\sum_{k,l}\frac{1}{\lmd}(BA)_{ik}(B^tA^t)_{lj}\dd\ot v_k^*\ot v_l-x\ot v^*_i\ot v_j,\\
	
	\mc{G}\ot(V_{A,B}^*\ot V_{A,B}) &\longrightarrow &\mc{G}\op \mc{G}\ot(V_{A,B}^*\ot V_{A,B})\\
	x\ot v_i^*\ot v_j&\longmapsto & (A^tB)_{ji}x+x\ot v_i^*\ot v_j+\sum_{k,l}(uB^t)_{li}(B^{-1})_{kj}\ot v_k^*\ot v_l,\\
	
\end{array}$$
$$\psi^t_4:\mc{G} \op \mc{G}\ot(V_{A,B}^*\ot V_{A,B})\longrightarrow \mc{G}$$
is defined by
$$\begin{array}{rcl}
	\mc{G} &\longrightarrow &\mc{G}\\
	x&\longmapsto &  x (\dd-1),\\
	
	\mc{G}\ot(V_{A,B}^*\ot V_{A,B}) &\longrightarrow &\mc{G}\\
	x\ot v_i^*\ot v_j&\longmapsto & x[(A^tB)_{ji}-(AuB^t)_{ji}].
\end{array}$$
Similar to the construction of the resolution (\ref{eq 4}), we have the following free resolution of $\kk$ as a left $\mc{G}$-module. 
$$0 \ra \mc{G}\xra{\vph_4} \mc{G}\ot (W_{A,B}^*\ot W_{A,B}) \op \mc{G}\xra{\vph_3} \mc{G}\ot(W_{A,B}^*\ot W_{A,B})\op \mc{G}\ot(V_{A,B}^*\ot V_{A,B})\hspace{10mm}$$
$$\hspace{60mm}\xra{\vph_2} \mc{G}\op \mc{G}\ot(V_{A,B}^*\ot V_{A,B})\xra{\vph_1}\mc{G}\ra 0,$$
where the morphisms are defined as follows,
$$\begin{array}{rcl}
	\vph_4:\mc{G}&\longrightarrow &\mc{G}\ot(W_{A,B}^*\ot W_{A,B})  \op \mc{G}\\
	x&\longmapsto &\sum_{i,j}  x[(AB^t)_{ji}-(A^tuB)_{ji}]\ot w_i^*\ot w_j+x(\dd-1),
\end{array}$$
$$\vph_3:  \mc{G}\ot(W_{A,B}^*\ot W_{A,B}) \op \mc{G} \longrightarrow  \mc{G}\ot(W_{A,B}^*\ot W_{A,B})\op \mc{G}\ot(V_{A,B}^*\ot V_{A,B})$$ is defined  by
$$\begin{array}{rcl}
	\mc{G}\ot(W_{A,B}^*\ot W_{A,B})  &\longrightarrow & \mc{G}\ot(W_{A,B}^*\ot W_{A,B})\op \mc{G}\ot(V_{A,B}^*\ot V_{A,B})\\
	x\ot w_i^*\ot w_j &\longmapsto &  x\ot w_i^*\ot w_j+x\sum_{k,l} (A^tu)_{li}(A^{-1})_{kj} \ot w_k^*\ot w_l\\
	&&\hspace{5mm}+x\sum_{k,l}\frac{1}{\lmd}(B^tA^t)_{ki}(BA)_{jl}\dd \ot v_k^*\ot v_l-  x\ot v_i^*\ot v_j,\\
	\mc{G} &\longrightarrow &\mc{G}\ot(W_{A,B}^*\ot W_{A,B})\op \mc{G}\ot(V_{A,B}^*\ot V_{A,B})\\
	x&\longmapsto & x(AB^t)_{ji}\ot w^*_i\ot w_j+ x ((A^tuB)_{ji}-(AB^t)_{ji})\ot v^*_i\ot v_j,\\
\end{array}$$
$$\vph_2: \mc{G}\ot(W_{A,B}^*\ot W_{A,B}) \op \mc{G}\ot(V_{A,B}^*\ot V_{A,B})\longrightarrow \mc{G}\op \mc{G}\ot(V_{A,B}^*\ot V_{A,B})$$ is defined by
$$\begin{array}{rcl}
	\mc{G}\ot(W_{A,B}^*\ot W_{A,B}) &\longrightarrow &\mc{G}\op \mc{G}\ot(V_{A,B}^*\ot V_{A,B})\\
	x\ot w_i^*\ot w_j&\longmapsto & x(u_{ji}-\dt_{ji})\\
	&&+x\sum_{k,l}\frac{1}{\lmd}(B^tA^t)_{ki}(BA)_{jl}\dd v_k^*\ot v_l-  x\ot v_i^*\ot v_j,\\
	
	\mc{G}\ot(V_{A,B}^*\ot V_{A,B}) &\longrightarrow &\mc{G}\op \mc{G}\ot(V_{A,B}^*\ot V_{A,B})\\
	x\ot v_i^*\ot v_j&\longmapsto & -x\dt_{ji}-x\ot v_i^*\ot v_j-x\sum_{k,l}(A^tu)_{li} A^{-1}_{kj}\ot v^*_k\ot v_l,\\
\end{array}$$
$$\vph_1:\mc{G} \op \mc{G}\ot(V_{A,B}^*\ot V_{A,B})\longrightarrow \mc{G}$$
is defined by
$$\begin{array}{rcl}
	\mc{G} &\longrightarrow &\mc{G}\\
	x&\longmapsto &  x (\dd-1),\\
	
	\mc{G}\ot(V_{A,B}^*\ot V_{A,B}) &\longrightarrow &\mc{G}\\
	x\ot v_i^*\ot v_j&\longmapsto & x(\dt_{ji}-u_{ji}).
\end{array}$$

We write $Q_0=Q_4=\mc{G}$, $Q_3=\mc{G}\ot (W_{A,B}^*\ot W_{A,B})\op \mc{G}$, $Q_2=\mc{G}\ot (W_{A,B}^*\ot W_{A,B})\op \mc{G}\ot (V_{A,B}^*\ot V_{A,B})$, $Q_1=\mc{G}\op\mc{G}\ot (V_{A,B}^*\ot V_{A,B})$. Let $\nu$ be the automorphism defined by $\nu(u)= A^{-1}A^t uB (B^t)^{-1}$ and  $\nu(\dd ^{\pm 1})=\dd ^{\pm 1}$. The following diagram is commutative
$$\xymatrix
{0 \ar[r]& \ar^{f_4}[d]Q_4\ar^-{\psi^t_1}[r]&\ar^{f_3}[d] Q_3\ar^-{\psi^t_2}[r]&\ar^{f_2}[d] Q_2\ar^-{\psi^t_3}[r] &\ar^{f_1}[d]Q_1\ar^-{\psi^t_4}[r]&\ar^{f_0}[d]Q_0\ar[r]&0\\
	0 \ar[r]& {}_\nu Q_4\ar^-{\vph_4}[r]& {}_\nu Q_3\ar^-{\vph_3}[r]& {}_\nu Q_2\ar^-{\vph_2}[r] &{}_\nu Q_1\ar^-{\vph_1}[r]&{}_\nu Q_0\ar[r]&0}$$
where the morphisms $f_i$, $i=0,1,\cdots, 4$,  are defined as follows. 
$$\begin{array}{rcl}
	f_4:\mc{G} &\longrightarrow & {}_\nu\mc{G}\\
	x &\longmapsto &  -\nu(x),\\
\end{array}$$
$$	f_3:\mc{G}\ot(W^*\ot W) \op \mc{G} \longrightarrow {}_\nu\mc{G}\ot(W^*\ot W) \op {}_\nu\mc{G}$$
is defined by
$$\begin{array}{rcl}
	\mc{G}\ot(W^*\ot W)  &\longrightarrow & {}_\nu\mc{G}\ot(W^*\ot W) \op {}_\nu\mc{G}\\
	x\ot w_i^*\ot w_j &\longmapsto & \sum_{p,q=1}^n -B_{pi}A_{qj}\nu(x)\ot w_p^*\ot w_q,\\
	\mc{G} &\longrightarrow & {}_\nu\mc{G}\ot(W^*\ot W) \op {}_\nu\mc{G}\\
	x &\longmapsto &  -\nu(x),\\
\end{array}$$
$$f_2:\mc{G}\ot(W^*\ot W)\op\mc{G}\ot(V^*\ot V) \longrightarrow  {}_\nu\mc{G}\ot(W^*\ot W) \op {}_\nu\mc{G}\ot(V^*\ot V)$$is defined by
$$\begin{array}{rcl}
	\mc{G}\ot(W^*\ot W) &\longrightarrow & {}_\nu\mc{G}\ot(W^*\ot W) \op {}_\nu\mc{G}\ot(V^*\ot V)\\
	x\ot w_i^*\ot w_j &\longmapsto & \sum_{p,q=1}^n B_{pi}A_{qj}\nu(x)\ot w_p^*\ot w_q,\\
	\mc{G}\ot(V^*\ot V)&\longrightarrow & {}_\nu\mc{G}\ot(W^*\ot W) \op {}_\nu\mc{G}\ot(V^*\ot V)\\
	x\ot v_i^*\ot v_j &\longmapsto & -\sum_{p,q=1}^n B_{pi}A_{qj}\nu(x)\ot v_p^*\ot v_q,\\
\end{array}$$
$$f_1: \mc{G}\op\mc{G}\ot(V^*\ot V) \longrightarrow  {}_\nu\mc{G}\op{}_\nu\mc{G}\ot(V^*\ot V)$$
is defined by
$$\begin{array}{rcl}
	\mc{G}\ot(V^*\ot V)  &\longrightarrow & {}_\nu\mc{G}\op{}_\nu\mc{G}\ot(V^*\ot V)\\
	x\ot v_i^*\ot v_j &\longmapsto & \sum_{p,q=1}^n B_{pi}A_{qj}\nu(x)\ot v_p^*\ot v_q,\\
	\mc{G} &\longrightarrow & {}_\nu\mc{G}\op{}_\nu\mc{G}\ot(V^*\ot V)\\
	x &\longmapsto &  \nu(x),\\
\end{array}$$
$$\begin{array}{rcl}
	f_0:\mc{G} &\longrightarrow & \op{}_\nu\mc{G}\\
	x &\longmapsto &  \nu(x).\\
\end{array}$$

Therefore, we obtain that $$\Ext_\mc{G}(\kk,\mc{G})\cong\begin{cases}0,& i\neq 4;
\\{}_\eta\kk,&i=4,\end{cases}$$
where $\eta$ is the algebra homomorphism $\vps\nu: \mc{G}\to \kk$, that is,  
\[\eta(u)=A^{-1}A^t B (B^t)^{-1}\quad \text{and}\quad \eta(\dd^{\pm 1})=1.\] 
So the algebra $\mc{G}$ is right AS-Gorentein. By  Lemma \ref{lem AS}, we obtain that the algebra $\mc{G}$ is a twisted CY algebra with a Nakayama automorphism $S^{-2}[\eta S]^r$. We have that $S^{-2}(u)=\dd (A^t)^{-1}AuA^{-1}A^t\dd^{-1}$ and $S^{-2}(\dd^{\pm 1})=\dd^{\pm 1}$. So $S^{-2}[\eta  S]^r$ is the  automorphism defined $(S^{-2}[\eta S]^r)(u)=\dd (A^t)^{-1}A u B^tB^{-1}\dd^{-1}$ and $(S^{-2}[\eta S]^r)(\dd^{\pm 1})=\dd^{\pm 1}$.  Let $\mu$ be the automorphism defined by $\mu(u)= (A^t)^{-1}A u B^tB^{-1}$ and $\mu(\dd^{\pm 1})=\dd^{\pm 1}$. It differs from  $S^{-2}[\eta S]^r$ only by an inner automorphism, so the algebra $\mu$ is also a Nakayama automorphism of $\mc{G}(A,B)$. \qed

We end this section by discussing the CY property of the Hopf-Galois objects. \cite[Theorem 5.3,5.4]{mro} gives the classification of the Galois objects over $\mc{G}(A,B)$.
Let $A,B\in \GL_n(\kk)\;(n\le 2)$ such that $B^tA^tBA=\lambda I_n$ for some $\lmd\in \kk^\times$. The left $\mc{G}(A,B)$-Galois objects are the algebras $\mc{G}(A,B|C,D)$ such that  $C,D\in \GL_m(\kk)$ ($m\le 2$) satisfying $D^tC^tDC=\lmd I_m$ and $\tr (AB^t)=\tr(CD^t)$.  

	 As a corollary of  Theorem \ref{thm main1}, we obtain the following result. 
\begin{thm}\label{thm galois}
	Let $m,n\in \NN , n,m\le 2$ and let $A,B\in \GL_n(\kk),C,D\in \GL_m(\kk)$ such that $B^tA^tBA=\lmd I_n$, $D^tC^tDC=\lmd I_m$ and $\tr (AB^t)=\tr(CD^t)=\mu$ for some $\lmd, \mu\in \kk^{\times}$. The algebra $\mc{G}(A,B|C,D)$ is twisted CY with a Nakayama automorphism $\mu$ defined by $\mu(u)= (A^t)^{-1}A u D^tD^{-1}$ and $\mu(\dd^{\pm 1})=\dd^{\pm 1}$.
\end{thm}
\proof It follows from  Theorem \ref{thm main1} and Lemma \ref{lem AS}  that the algebra $\mc{G}(A,B)$ is AS-Gorenstein with left homological integral $\xi$ defined by $\xi(u)= (A^t)^{-1}A B^tB^{-1}$ and $\xi(\dd^{\pm 1})=1$.  Following from \cite[Theorem 2.5]{yu}, we conclude  that the algebra $\mc{G}(A,B|C,D)$ is a twisted CY algebra with a Nakayama automorphism $\mu'=S_{CD,AB}S_{AB,CD}[\xi]^l_{AB,CD}$, where $[\xi]^l_{AB,CD}$ is the left winding automorphism of  $\mc{G}(A,B|C,D)$ defined by $[\xi]^l_{AB,CD}(x)=\xi(x^{AB,AB}_{(1)})x^{AB,CD}_{(2)}$, for any $x\in \mc{G}(A,B|C,D)$. In $\mc{G}(A,B|C,D)$, we have $$S_{CD,AB}S_{AB,CD}(u)=A^{-1}(B^t)^{-1}uD^tC,\quad S_{CD,AB}S_{AB,CD}(\dd^{\pm 1})=\dd^{\pm 1}.$$  
Therefore, the automorphism $\mu'=S_{CD,AB}S_{AB,CD}[\xi]^l_{AB,CD}$ is defined by $$\mu'(u)=(A^t)^{-1}B^{-1}uD^tC,\quad  \mu'(\dd^{\pm 1})=\dd^{\pm 1}.$$
 Let $\mu$ be the algebra automorphism defined by $$\mu(u)= (A^t)^{-1}A u D^tD^{-1},\quad \mu(\dd^{\pm 1})=\dd^{\pm 1}.$$ In $\mc{G}(A,B|C,D)$, we also have $\dd^{-1}u\dd=BAuC^{-1}D^{-1}$. Therefore, the algebra automorphisms $\mu$ and $\mu'$ are differ only by an inner automorphism. So the algebra automorphism $\mu$ is also a Nakayama automorphism of $\mc{G}(A,B|C,D)$. 

\section{Gerstenhaber-Schack cohomology}
Gerstenhaber-Schack cohomology of a Hopf algebra with coefficients in Hopf bimodules was introduced in \cite{gs1, gs2} by using an explicit bicomplex. Since the category of  Hopf bimodules and the category of Yetter-Drinfeld modules are equivalent, we can work with Yetter-Drinfeld modules. Let $H$ be a Hopf and $V$ a Yetter-Drinfeld module over $A$. The Gerstenhaber-Schack cohomology of $A$ with coefficients in $V$ is usually denoted by $\H^*_{GS}(H,V)$. A special instance of Gerstenhaber-Schack cohomology is bialgbra cohomology, given by $\H_b^*(H)=\H^*_{GS}(H,\kk)$. 

Taillefer proved that Gerstenhaber-Schack cohomology is in fact an $\Ext$-functor (\cite{ta1, ta2})
\begin{equation}\label{eq 5}\H^*_{GS}(H,V)\cong \Ext_{\mc{YD}^H_H}^*(\kk,V).\end{equation}
This description can be viewed as a defition of Gerstenhaber-Schack cohomology. 

The Gerstenhaber-Schack cohomological dimension of a Hopf algebra $H$ is defined to be 
$$\cd_{GS}={\rm sup}\{n:\H^*_{GS}(H,V)\neq0 \text{ for some }V\in \mc{YD}^H_H\}\in\NN \cup \{\infty\}$$

\begin{thm} Assume that $char(\kk)=0$. Let $A,B\in \GL_n(\kk)\;(n\le 2)$ such that $B^tA^tBA=\lambda I_n$ for some $\lmd\in \kk^\times$ and such that any solution of the equation $X^2-\sqrt{\lmd^{-1}}\tr(AB^t)X+1=0$ is generic. 	Then we can obtain the followings:
\begin{enumerate}
	\item[(i)] The bialgebra cohomology of $\mc{G}(A,B)$ is as follows	$$\H_b^{n}(\mc{G}(A,B))\cong \begin{cases}
		0&n\neq 0,1,3,4\\
		\kk &n=0,1,3,4.
	\end{cases}$$
\item[(ii)] The Gerstenhaber-Schack cohomological dimension of $\mc{G}(A,B)$ is 4.
\end{enumerate}
	\end{thm}
\proof  Let $\mc{G}=\mc{G}(A,B)$. 

(i) We have  $$\H_b^*(\mc{G})=\H^*_{GS}(\mc{G},\kk)\cong\Ext_{\mc{YD}^{\mc{G}}_{\mc{G}}}^*(\kk,\kk).$$
By the assumption, the algebra $\mc{G}(A,B)$ is cosemisimple. So every object in $\mc{M}^{\mc{G}}$ is projective. By \cite[Proposition 3.1 and Corollary 3.4]{bi}, every free Yetter-Drinfeld module is a projective object in $\mc{YD}^{\mc{G}}_{\mc{G}}$ and the category  $\mc{YD}^{\mc{G}}_{\mc{G}}$ has enough projective objects. Therefore, the bialgebra cohomology of $\mc{G}$ is  the cohomology of the complex obtained by applying the functor $\Hom_{\mc{YD}^{\mc{G}}_{\mc{G}}}(-,\kk)$ to any projective resolution of the trivial Yetter-Drinfeld module $\kk$. We can use the  resolution (\ref{eq 4}) in Theorem \ref{main thm}. We keep the notation as there. Let ${\bf P}_.$ denote the following complex
$$0 \!\ra\! \kk\boxtimes\mc{G}\!\xra{\psi_4}\! (V_{A,B}^*\ot V_{A,B})\boxtimes \mc{G} \op \kk\boxtimes\mc{G}\!\xra{\psi_3} \!(V_{A,B}^*\ot V_{A,B})\boxtimes \mc{G}\op (W_{A,B}^*\ot W_{A,B})\boxtimes \mc{G}$$
$$\hspace{50mm}\xra{\psi_2} \kk\boxtimes\mc{G}\op(W_{A,B}^*\ot W_{A,B})\boxtimes \mc{G}\xra{\psi_1}\kk\boxtimes\mc{G} \ra 0$$
It is easy to see that $$\Hom_{\mc{YD}^{\mc{G}}_{\mc{G}}}(\kk\boxtimes\mc{G},\kk)\cong \Hom_{\mc{M}^{\mc{G}}}(\kk,\kk)\cong \kk.$$
Both $V_{A,B}$ and $W_{A,B}$ are the fundamental comodules. We have
$$\Hom_{\mc{YD}^{\mc{G}}_{\mc{G}}}((V_{A,B}^*\ot V_{A,B})\boxtimes \mc{G},\kk)\cong \Hom_{\mc{M}^{\mc{G}}}(V_{A,B}^*\ot V_{A,B},\kk)\cong \kk.$$
Therefore, $\Hom_{\mc{YD}^{\mc{G}}_{\mc{G}}}({\bf P}_.,\kk)$ is isomorphic to the following complex
$$0\ra\kk\xra{\left(\begin{array}{cc}0&0\end{array}\right)} \kk^2\xra{\left(\begin{array}{cc}1&0\\1&0\end{array}\right)}\kk^2\xra{\left(\begin{array}{cc}0&0\\1&1\end{array}\right)} \kk^2\xra{\left(\begin{array}{c}0\\0\end{array}\right)} \kk\ra0.$$
By compute the cohomology, we obtain the announced result for the bialgebra cohomology of $\mc{G}(A,B)$.

(ii) The algebra $\mc{G}(A,B)$ is cosemisimple, as mentioned in the proof of (i), every free Yetter-Drinfeld module is a projective object in $\mc{YD}^{\mc{G}}_{\mc{G}}$ and the category  $\mc{YD}^{\mc{G}}_{\mc{G}}$ has enough projective objects. Let $V$ be a Yetter-Drinfeld module, the Gerstenhaber-Schack cohomology $\H^*_{GS}(H,V)$ is the cohomology of the complex obtained by applying the functor $\Hom_{\mc{YD}^{\mc{G}}_{\mc{G}}}(-,V)$ to any projective resolution of the trivial Yetter-Drinfeld module $\kk$. By Theorem \ref{main thm}, the trivial module $\kk$ has a projective resolution of length 4. So $\cd_{GS}\mc{G}\se 4$. However, bialgebra cohomology computed in (i) grantees that $\cd_{GS}\mc{G}\le 4$. In conclusion, we obtain that $\cd_{GS}\mc{G}=4$. \qed

\subsection*{Acknowledgement}  Yu is supported by grant NSFC (No. 11871186). Wang was partially supported by Simons collaboration grant \#688403 and Air Force Office of Scientific Research grant FA9550-22-1-0272. 

\newpage

\vspace{5mm}

\bibliography{}

\end{document}